\documentclass[preprint,review,12pt]{elsarticle}

\usepackage[utf8]{inputenc} 
\usepackage[T1]{fontenc}    
\usepackage{hyperref}       
\usepackage{url}            
\usepackage{booktabs}       
\usepackage{amsfonts}       
\usepackage{amsmath, amssymb}
\usepackage{nicefrac}       
\usepackage{microtype}      
\usepackage{lipsum}
\usepackage{natbib}         
\usepackage{tabularx}
\usepackage{multirow}       
\usepackage{graphicx}
\usepackage{caption}
\usepackage{subcaption}
\usepackage{makecell}       
\usepackage{float}          
\usepackage{gensymb}        
\usepackage{epsfig}
\usepackage[flushleft]{threeparttable}
\usepackage{lineno}         
\usepackage{textgreek}      
\usepackage[dvipsnames]{xcolor}
\usepackage{algorithm}
\usepackage{algpseudocode}
\usepackage{bbm}            
\graphicspath{ {./Figs/} }

\newcommand{\decreases}{\mathbin{\text{\rotatebox[origin=c]{180}{$\simeq$}}}}

\usepackage{geometry}
\begin{document}

\begin{frontmatter}

\title{A probabilistic peridynamic framework with an application to the study of the statistical size effect}

\author[Exeter,Turing]{Mark Hobbs\corref{mycorrespondingauthor}}
\cortext[mycorrespondingauthor]{Corresponding author}
\ead{mhobbs@turing.ac.uk}

\author[Exeter]{Hussein Rappel}
\ead{h.rappel@exeter.ac.uk}

\author[Exeter,digiLab]{Tim Dodwell}
\ead{tim@digilab.co.uk}

\address[Exeter]{Department of Engineering, Faculty of Environment, Science and Economy, University of Exeter, United Kingdom}
\address[Turing]{The Alan Turing Institute, British Library, 96 Euston Road, London NW1 2DB, United Kingdom}
\address[digiLab]{digiLab, Kings Wharf, The Quay, Exeter, EX2 4AN}


\begin{abstract}

Mathematical models are essential for understanding and making predictions about systems arising in nature and engineering. Yet, mathematical models are a simplification of true phenomena, thus making predictions subject to uncertainty. Hence, the ability to quantify uncertainties is essential to any modelling framework, enabling the user to assess the importance of certain parameters on quantities of interest and have control over the quality of the model output by providing a rigorous understanding of uncertainty. Peridynamic models are a particular class of mathematical models that have proven to be remarkably accurate and robust for a large class of material failure problems. However, the high computational expense of peridynamic models remains a major limitation, hindering \textit{`outer-loop'} applications that require a large number of simulations, for example, uncertainty quantification. This contribution provides a framework to make such computations feasible. By employing a Multilevel Monte Carlo (MLMC) framework, where the majority of simulations are performed using a coarse mesh, and performing relatively few simulations using a fine mesh, a significant reduction in computational cost can be realised, and statistics of structural failure can be estimated. The results show a speed-up factor of 16$\times$ over a standard Monte Carlo estimator, enabling the forward propagation of uncertain parameters in a computationally expensive peridynamic model. Furthermore, the multilevel method provides an estimate of both the discretisation error and sampling error, thus improving the confidence in numerical predictions. The performance of the approach is demonstrated through an examination of the statistical size effect in quasi-brittle materials.

\end{abstract}


\end{frontmatter}


\section{Introduction}

Design approaches remain broadly the same across different engineering disciplines (e.g. aerospace, structural and mechanical). All disciplines depend heavily upon empirical formulas and large safety factors. This approach leads to highly conservative designs with low material utilisation. The benefits of improving material utilisation are clear (e.g. lighter vehicles achieve greater fuel efficiency) but certifying the safety and reliability of novel structural forms requires expensive programmes of testing. As the demand for more efficient structures increases, the need for new design approaches becomes more pressing. When experimental data is incomplete, a better approach to examine structural reliability might be provided through numerical simulations and stochastic methods. 

Uncertainties in structural analysis arise from multiple sources, for example, material properties, loading conditions and geometry \cite{Schueller2007}. Current design approaches generally rely on deterministic models, and large safety factors must be applied to account for the inherent uncertainties. Uncertainties can be examined by using stochastic simulation methods, where uncertain input parameters, such as material properties, are treated as random variables. Methods for the forward propagation of uncertainty (where sources of uncertainty are propagated through a model to evaluate the uncertainty in the output) can be broadly classified as intrusive or non-intrusive \cite{Lin2012}. Intrusive uncertainty quantification (UQ) methods reformulate the original deterministic governing equations that describe the physical process. Non-intrusive UQ methods sample uncertain input parameters from a probability distribution and the deterministic governing equations are solved for each sample. The output is a distribution of the quantity of interest (QoI) from which various statistics, such as the mean and variance, can be computed. Non-intrusive methods are viewed in terms of inputs and outputs and the form of the governing equations is irrelevant. This makes non-intrusive methods very general.

Monte Carlo (MC) methods are a broad class of approaches that repeatedly sample random input variables to approximate the statistics of the response of a model. Monte Carlo simulation remains the most popular and straightforward method to quantify the effect of uncertainty in systems with random input parameters. Random variables are sampled from a probability distribution, a deterministic model is solved for each sample, and the output is a distribution of the quantity of interest, for example, the maximum stress at failure. The popularity of Monte Carlo simulations can be explained by a number of factors: (1) MC simulations are non-intrusive and simple to implement, (2) MC simulations can be applied to any problem with a stochastic element (including non-linear and discontinuous problems), and (3) MC simulations can be applied to problems with non-Gaussian random variables. 

The order of convergence of MC simulations is slow and computing the statistics of the QoI requires a large number of model evaluations. If we are interested in calculating the expected value $\mathbb{E}[Q]$ of a quantity of interest $Q$, by the central limit theorem, the rate of convergence of MC is $\mathcal{O}(1 / \sqrt{N})$, where $N$ is the number of samples. If the computational cost of a single deterministic simulation is high, standard Monte Carlo simulation is prohibitively expensive. This limitation has led to the development of different strategies to accelerate Monte Carlo simulations. These strategies, all based on the idea of variance reduction, can be broadly classified into three categories: (1) techniques that reduce the total number of required samples by carefully selecting the samples to reduce the error, i.e. they are not chosen randomly and independently. These techniques are generally known as quasi-Monte Carlo (QMC) methods and can achieve a convergence rate close to $\mathcal{O}(1 / N)$. For a review of QMC methods, the reader is referred to \citet{Dick2013}. (2) control variate strategies that exploit information about the error in estimates of random variables with a known expected value to reduce the error of an estimate of an unknown random variable. Control variate methods are useful when a simple version of the problem is available and can be easily solved. (3) Multilevel methods that utilise a low accuracy but computationally cheap model for the majority or samples, and employ a high accuracy and computationally expensive model for a small subset of samples. 

In this work, we employ the multilevel Monte Carlo (MLMC) method. The aim of MLMC is attain the same solution error as MC but at a significantly reduced computational cost. The standard MC estimator is computationally expensive as all samples must be computed using a fine mesh that guarantees a small discretisation error. A significant reduction in computational cost can be realised by taking the majority of samples on a coarse mesh (low accuracy but computationally cheap), and taking relatively few samples on a fine mesh (high accuracy but computationally expensive). This is made possible by isolating the error sources in the estimator: (1) sampling error (variance) and (2) discretisation error (deterministic error). The sampling error is controlled by using a low accuracy but computationally cheap model to take a large number of samples, and the discretisation error is reduced to a defined tolerance by employing a sufficiently fine mesh.

Multilevel techniques were first introduced by \citet{Heinrich1999} and \citet{Heinrich2001} and later popularised by \citet{Giles2008} for option pricing in computational finance. \citet{Cliffe2011} were the first to apply multilevel methods in the field of engineering, motivated by the study of uncertainty in groundwater flows. Since \citet{Cliffe2011} recognised the potential of multilevel methods, there has been a wide range of applications in engineering and scientific fields, for example, \citet{Dodwell2021} employed MLMC to estimate the probability of failure of composite materials, and \citet{Clare2022} assessed the risk of coastal flooding. For a detailed review of multilevel Monte Carlo methods, the reader is referred to the work of \citet{Giles2015}.

The peridynamic theory of solid mechanics, introduced by \citet{Silling2000}, is an integral-type non-local theory of solid mechanics that provides a robust theoretical framework for developing numerical models capable of simulating the failure behaviour of a wide range of materials. The peridynamic model defines material behaviour at a point in a continuum body as an integral equation of the surrounding displacement. This is in contrast to the classical theory of solid mechanics, where the material behaviour at a point is defined by partial differential equations. The classical theory is only valid if the body under analysis has a spatially continuous and differentiable displacement field. Spatial derivatives are not defined across discontinuities, and the classical theory breaks down when a body fractures. The governing equations of peridynamics do not require a spatially continuous and differentiable displacement field and the initiation and propagation of damage is an emergent behaviour. No additional assumptions or techniques are required for modelling damage and fracture. The high computational cost of peridynamic simulations remains a major limitation, and as a consequence, the application of forward UQ methods that require a large number of repeat runs is prohibitively expensive. To the best of the authors knowledge, there is no work within the peridynamic literature that examines the forward propagation of uncertain parameters. This work demonstrates that forward uncertainty quantification for expensive peridynamic simulations is feasible by employing the multilevel Monte Carlo method. 


The aims of this paper are: (1) demonstrate that forward uncertainty quantification for expensive peridynamic simulations is feasible by employing the MLMC method, (2) quantify the computational savings, (3) demonstrate the importance of forward UQ by selecting examples where a deeper understanding of the physical behaviour can be gained by considering uncertainty. We focus on quasi-brittle materials because the range of experimental data is greater than that of any other material and quasi-brittle materials exhibit a significant size effect. A stochastic model is required for a complete examination of the mechanisms that govern the structural size effect. This work also provides new insights into the convergence behaviour of bond-based peridynamic models. To the best of the authors knowledge, convergence studies of the predicted structural response are missing from the peridynamics literature. Existing convergence studies only consider static elastic problems \cite{Seleson2016}. 

The paper is organised as follows: Section \ref{section:peridynamic_theory} introduces the peridynamic theory and Section \ref{section:numerical_model} briefly describes the numerical model (a bond-based peridynamic model). Section \ref{section:mlmc} details the standard and multilevel Monte Carlo methodology. Section \ref{section:random_fields} describes the modelling of material properties as spatially correlated random fields. Section \ref{section:case_studies} presents two case studies. The presented problems have been selected as examples where a deeper understanding of the physical behaviour can be gained by considering uncertainty. Section \ref{section:cdf} demonstrates a method for computing the cumulative distribution function using samples obtained on multiple mesh levels. Section \ref{section:discussion} discusses the results and Section \ref{section:conclusions} concludes the paper. 


\section{Peridynamic theory}\label{section:peridynamic_theory}

The peridynamic theory, introduced by \citet{Silling2000} in 2000, is a non-local theory of solid mechanics that is formulated in terms of integral equations rather than partial differential equations. The governing equations do not require a spatially continuous and differentiable displacement field and damage localisation and fracture naturally emerge. No additional assumptions or techniques are required for modelling damage and fracture.

There are two primary formulations of the peridynamic theory: bond-based \cite{Silling2000} and state based theory \cite{Silling2007}. In the original bond-based theory, the Poisson's ratio is limited to a fixed value. \citet{Silling2007} later introduced a generalised state-based theory that overcomes the limitations of the original theory. This paper employs the bond-based theory due to its lower computational expense and proven predictive capabilities. 

\subsection{Peridynamic continuum model}

The bond-based peridynamic theory is briefly presented. It is not the purpose of this work to explain the peridynamic theory in detail and the reader is referred to \cite{Silling2000,Silling2005,Silling2010} for a comprehensive treatment of the theory. A mechanically intuitive but less rigorous way of obtaining the governing equations can be found in \cite{Bobaru2017}.

\begin{figure}[H]
    \centering
    \includegraphics[width = 0.75\textwidth]{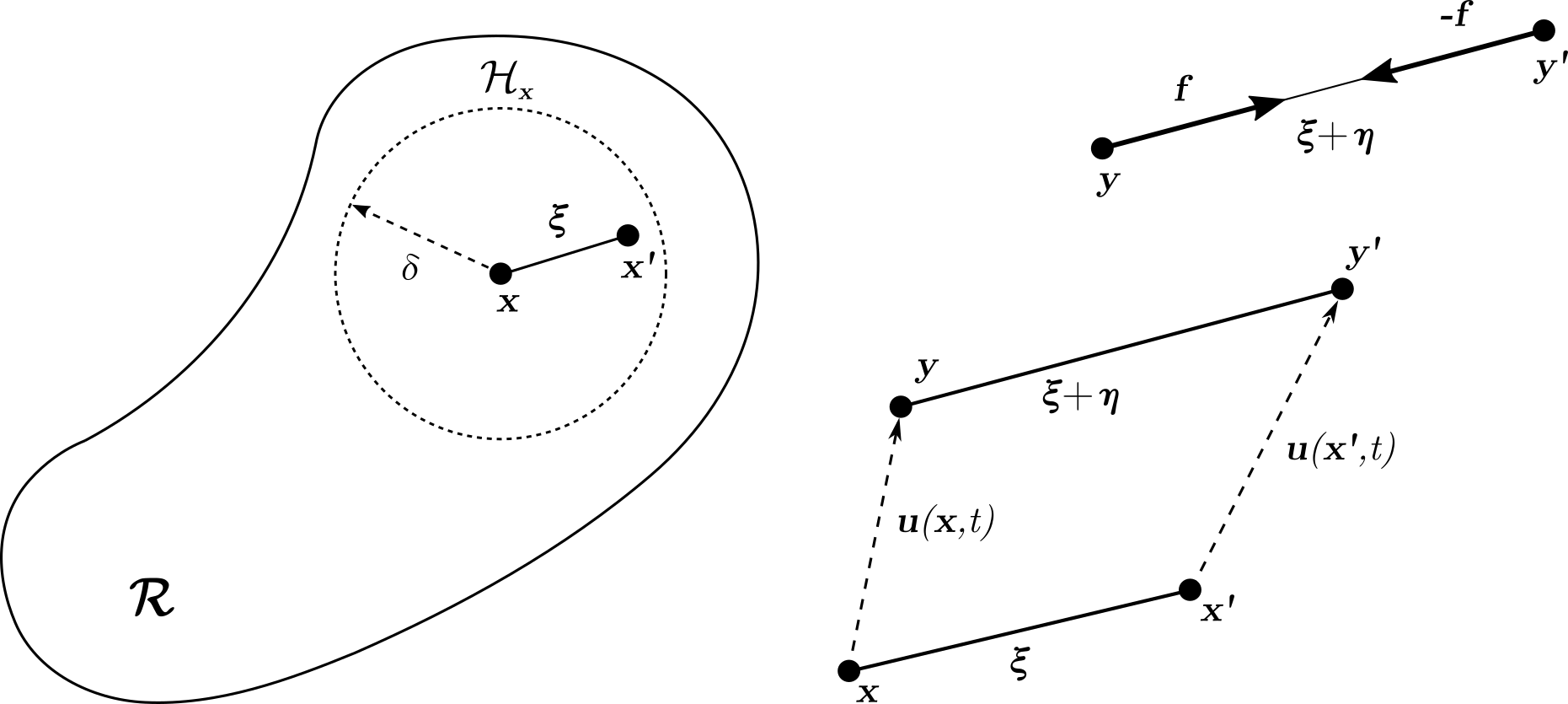}
    \caption{Peridynamic continuum and kinematics of particle pair and bond-based pairwise force function}
    \label{fig:peridynamic_continuum}
\end{figure}

Assuming that a body occupies a spatial region \(\mathcal{R}\), for any material point \(\textbf{x} \in \mathcal{R}\), a pairwise force function \(\textbf{f}\) can be defined to describe the interaction between particles within a finite distance \(\delta\) of \(\textbf{x}\), at any time \(t\), where \(\textbf{u}\) represents the displacement of a material point (see Fig. \ref{fig:peridynamic_continuum}).

\begin{equation}\label{eq:PD_Pairwise_Force_Function}
    \textbf{f} = \textbf{f}(\textbf{x},\textbf{x}',\textbf{u}(\textbf{x},t),\textbf{u}(\textbf{x}',t),t), \quad \textbf{x}' \in \mathcal{R}: || \textbf{x}' - \textbf{x} || \leqslant \delta
\end{equation}

The peridynamic equation of motion for a single material point \(\textbf{x}\) at a point in time \(t\) is given by Newton's second law of motion, and is defined by Eq. (\ref{eq:PD_EOM_Dynamic}).

\begin{equation}\label{eq:PD_EOM_Dynamic}
    \rho\ddot{\textbf{u}}(\textbf{x},t) = \int_{\mathcal{H}_\textbf{x}} \textbf{f}(\textbf{u}(\textbf{x}',t) - \textbf{u}(\textbf{x},t), \textbf{x}' - \textbf{x})dV_{\textbf{x}'} + \textbf{b}(\textbf{x},t)
\end{equation}

\noindent $\rho$ is mass density, $\ddot{\textbf{u}}$ is particle acceleration, $\textbf{b}$ is body force per unit volume, and $\mathcal{H}_\textbf{x}$ is the neighbourhood of material point $\textbf{x}$. The size of the neighbourhood is defined by the horizon length $\delta$. For a 3D problem, the material point neighbourhood will be a sphere, and for a 2D problem, the neighbourhood will be circular.

\begin{equation}\label{eq:PD_Horizon}
    \mathcal{H}_\textbf{x} = \mathcal{H}_\textbf{x}(\textbf{x},\delta) = \{ \textbf{x}' \in \mathcal{R}: || \textbf{x}' - \textbf{x} || \leqslant \delta \}
\end{equation}

The pairwise force function \(\textbf{f}\) represents the force that particle \(\textbf{x}'\) exerts on particle \(\textbf{x}\) and contains all the constitutive information of the material under analysis. This interaction is commonly referred to as the peridynamic bond force. Particles separated by a distance greater than the horizon length \(\delta\) do not interact. The pairwise force function \(\textbf{f}\) is defined by Eq. (\ref{eq:PD_Force_Function}).

\begin{equation}\label{eq:PD_Force_Function}
    \textbf{f}(\boldsymbol{\eta},\boldsymbol{\xi}) = f(|\boldsymbol{\xi} + \boldsymbol{\eta}|,\boldsymbol{\xi}) \frac{\boldsymbol{\xi} + \boldsymbol{\eta}}{|\boldsymbol{\xi} + \boldsymbol{\eta}|}
\end{equation}

In a bond-based model, the force vector \(\textbf{f}\) is parallel to the deformed bond and the scalar bond force \(f\) (vector magnitude) is proportional to the bond stretch \(s\). The initial relative position vector of a pair of particles is denoted by \(\boldsymbol{\xi} = \textbf{x}' - \textbf{x}\), and the relative displacement vector is denoted by \(\boldsymbol{\eta} = \textbf{u}' - \textbf{u}\). The current relative position vector is given by \(\boldsymbol{\xi} + \boldsymbol{\eta}\).

To make a distinction between the peridynamic theory and other non-local theories, note that most non-local theories average some measure of strain within a neighbourhood of a material particle. The peridynamic theory dispenses with the concept of strain, which by its definition, requires the evaluation of partial derivatives of displacement \cite{Silling2005}.

\subsection{Non-locality}

The peridynamic theory is a non-local theory in which material points interact with each other directly over finite distances. This is in contrast to the classical theory of solid mechanics, where it is assumed that all forces are contact forces that act across zero distance (local theory). Physical justification of non-locality was provided by \citet{Bazant1991a}, and further discussion on the origins of non-locality (with a focus on the peridynamic theory), can be found in Chapter 1 of \citet{Bobaru2017} and \citet{Hobbs2021a}.

At the macroscale, the peridynamic horizon $\delta$ is a numerical constant with no physical meaning. This differentiates the peridynamic model from many numerical approaches, and the use of an ambiguous characteristic length parameter is avoided. For a given value of $\delta$, the parameters in a peridynamic model can be chosen to match a given set of physically measurable material properties. Therefore, an optimum value of $\delta$ must be chosen that provides high accuracy whilst balancing computational expense. Section \ref{section:numerical_convergence} discusses the selection of an optimum value of $\delta$. 

The reader should note the distinction between the non-local length scale in the peridynamic model (horizon $\delta$), and the non-local length scale in a spatially correlated random field (correlation length $l_c$). The correlation length $l_c$ is generally considered to be a material parameter reflecting the internal length scale of the microstructure. This will be discussed throughout the paper.


\section{Numerical model}\label{section:numerical_model}

To illustrate the framework, this work employs a two-dimensional bond-based peridynamic model. The aim of this work is to demonstrate the multilevel framework and a detailed treatment of the numerical model is not provided. All the results presented in this paper were obtained using an explicit scheme (as outlined in Fig. 4.14 of \citet{Hobbs2021a}). The reader is referred to \citet{Hobbs2021a} for implementation details. The main distinction of the model used in this work is the existence of two length scales: (1) the peridynamic horizon $\delta$ and (2) the correlation length $l_c$ in the random field. The generation of spatially correlated random fields is discussed in Section \ref{section:random_fields}. 

\subsection{Constitutive model}

It is generally assumed that the force-stretch ($f\textrm{-}s$) relationship of a peridynamic bond should be consistent with the macroscopic material response, and a failure mechanism is introduced into the model by eliminating the interaction between particle pairs when the stretch of the connecting bond exceeds a critical value. The stress-strain response of quasi-brittle materials is characterised by strain softening behaviour in the post-peak stage, and hence we employ the non-linear softening law, illustrated in Fig. \ref{fig:Decaying_exponential}, first proposed by \citet{Hobbs2021a}. The derivation of the model parameters for the two-dimensional plane stress and plane strain case are provided in Appendix A.

\subsection{Numerical convergence}\label{section:numerical_convergence}

\begin{figure}[H]
    \centering
    \includegraphics[width = 0.75\textwidth]{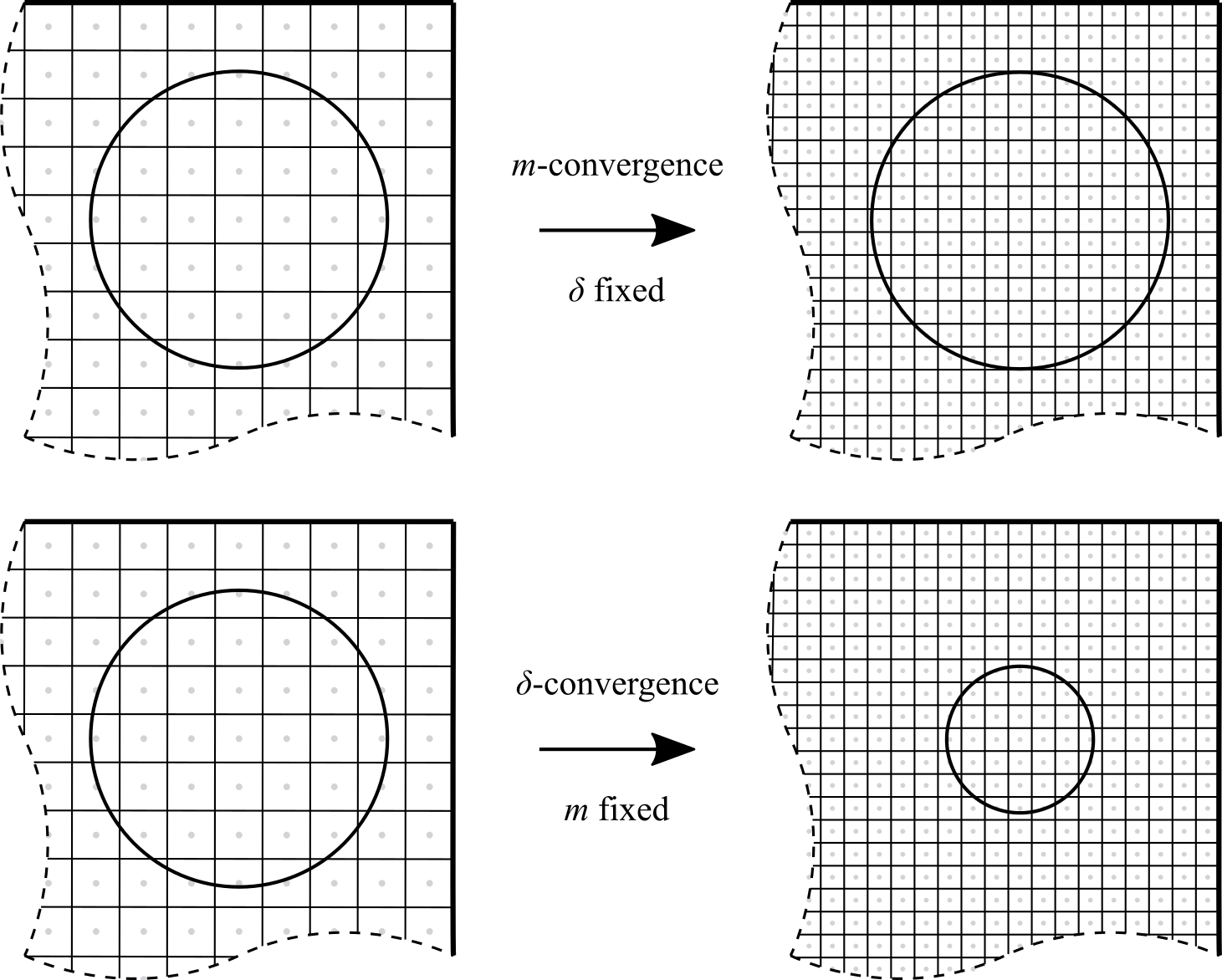}
    \caption[\(m\)-convergence and \(\delta\)-convergence]{Graphical representation of the two fundamental types of convergence: \(m\)-convergence and \(\delta\)-convergence.}
    \label{fig:convergence}
\end{figure}

The accuracy and convergence behaviour of a peridynamic model is complicated by the presence of a length scale. To determine an optimum value of \(\delta\), an additional parameter \(m\) must be introduced. \(m\) is the ratio between the horizon radius and grid resolution (\(m = \delta/\Delta x\)). \citet{Bobaru2009} and \citet{Ha2010} define and discuss two fundamental types of convergence: (1) \(m\)-convergence: \(\delta\) is fixed and \(m \rightarrow \infty\). This can also be stated as \(\delta\) is fixed and \(\Delta x \rightarrow 0\). (2) \(\delta\)-convergence: \(m\) is fixed and \(\delta \rightarrow 0\). This can also be stated as \(m\) is fixed and \(\Delta x \rightarrow 0\). See Fig. \ref{fig:convergence} for a graphical representation of the types of convergence. A third type of convergence can be defined: \(\delta m\)-convergence. This is a combination of \(\delta\)- and \(m\)-convergence. See \citet{Bobaru2009} for details.

In this work, we consider $\delta$-convergence, as it is generally agreed that $m$ should be close to 3. \citet{Madenci2014} investigated the choice of $m$ for macroscale problems and it was found that values of $m = 1$ and $m = 3$ achieved the highest accuracy when compared to the classical analytical solution for the displacement of a one-dimensional bar subjected to a defined initial strain. Values of $m$ much larger than 3 lead to excessive wave dispersion and become extremely computationally expensive. When fracture behaviour is also considered, values of $m$ less than 3 lead to grid dependence on crack propagation \cite{Bobaru2012, Madenci2014}. \citet{Hu2010} and \citet{Seleson2014} examined the $m$-convergence behaviour for two-dimensional models and \citet{Hobbs2021a} examined the $m$-convergence behaviour for three-dimensional models. Higher values of $m$ improve the spatial integration accuracy but $m \approx 3$ provides an acceptable approximation. A value of $m = \pi \; (\delta = \pi \Delta x)$ is generally recommended for macroscale problems and is found extensively throughout the literature. The $m$-ratio is set to $\pi$ for all problems in this paper.

\section{Multilevel Monte Carlo methodology}\label{section:mlmc}

The simplicity and dimensionality independence of MC methods makes them the most popular and straightforward technique to quantify the effect of uncertainty in systems with random input parameters. Principally, MC techniques are based on iterative calculations of a mathematical model for random samples of the parameters that describe the model. However, the iterative nature of these methods makes them inefficient for problems in which the model evaluations are computationally expensive, such as peridynamic models. Multilevel techniques are designed to overcome the disadvantages of conventional MC techniques by performing most simulations using a low accuracy but computationally cheap model, and relatively few simulations are performed using a high accuracy but computationally expensive model.

To explain the standard and multilevel Monte Carlo method, let us assume that we have a numerical model of a brittle/quasi-brittle structure that is subject to some uncertainty in the material properties. The accuracy and computational cost of the model is proportional to the number of degrees of freedom (\(M\)) and thus to the resolution of the mesh. Generally, we are interested in some scalar quantity of interest $Q = Q(\textbf{x}, t, \boldsymbol{\omega})$, for example, the failure load or the displacement at a particular point, where $\textbf{x}$ and $t$ are the spatial and temporal coordinates and $\boldsymbol{\omega}$ represents a vector of random variables that takes values in $\mathbb{R}^M$. $\boldsymbol{\omega}$ represents sources of uncertainty in the problem, in this case, the material properties. Note that the quantity of interest ($Q$) could be a function, for instance, the load-deflection response of a structure. For the presented case studies, the quantity of interest is the failure load, and the objective is to compute the expected value of $Q$, denoted $\mathbb{E}[Q]$, with a quantified level of uncertainty. However, for many real world applications, the probability distribution of $Q$ is of more interest. Methods for obtaining the probability distribution of $Q$ will be discussed. 

\subsection{Standard Monte Carlo simulation}

In a standard Monte Carlo (MC) simulation, a large number ($N$) of independent random realisations (or samples) of the parameters are generated. For every sample, the solution is computed using a numerical solver (finite element model, particle model etc). The accuracy of the solution is directly proportional to the resolution of the discretisation, and it is assumed that $Q_M \rightarrow Q$ as $M \rightarrow \infty$, therefore $\mathbb{E}[Q_M] \rightarrow \mathbb{E}[Q]$ for $M \rightarrow \infty$. The required accuracy and computational budget govern the selection of $M$. The standard Monte Carlo estimator for the expected value $\mathbb{E}[Q_M]$ of $Q_M$, based on $N$ samples, is given by Eq. (\ref{eq:MC_estimator}), where $Q_M^{(j)}$ is the quantity of interest of the $j^{th}$ sample.

\begin{equation}\label{eq:MC_estimator}
    \hat{Q}_{M, N}^{MC} = \frac{1}{N}\sum_{j=1}^{N} Q_M^{(j)}
\end{equation}

Note that $\hat{Q}_{M, N}^{MC}$ is an unbiased estimator of $\mathbb{E}[Q_M]$, meaning that $\mathbb{E}[\hat{Q}_{M, N}^{MC} ] = \mathbb{E}[Q_M]$. An estimator is \textit{unbiased} if its expectation is the quantify of interest that we wish to estimate. The accuracy of the estimator ($\hat{Q}_{M, N}^{MC}$) can be quantified using the \textit{root mean square error} (RMSE):

\begin{equation}\label{eq:MC_error}
    e(\hat{Q}_{M, N}^{MC}) = \sqrt{\mathbb{E}[(\hat{Q}_{M, N}^{MC} - \mathbb{E}[Q])^2]}
\end{equation}

An advantage of quantifying the accuracy of the estimator in this way is that the mean square error can be expanded and two distinct sources of error can be isolated: (1) the bias error and (2) the sampling error. 

\begin{equation}\label{eq:MC_error_expanded}
    e(\hat{Q}_{M, N}^{MC})^2 = \underbrace{\left(\mathbb{E}[Q_M - Q]\right)^2}_{bias \; error} + \underbrace{\frac{\mathbb{V}[Q_M]}{N}}_{sampling \; error}
\end{equation}

The first term in Eq. (\ref{eq:MC_error_expanded}) is the \textit{bias error} (sometimes referred to as the discretisation or numerical error). This arises as we are actually interested in the expected value $\mathbb{E}[Q]$ of $Q$, the unobtainable random variable corresponding to the exact solution without any numerical error. If we assume that the numerical model converges to the exact solution as the mesh is refined, $\mathbb{E}[Q_M] \rightarrow \mathbb{E}[Q]$ as $M \rightarrow \infty$, then we can state the following

\begin{equation}\label{eq:alpha_convergence}
    \mathbb{E}[Q_M - Q] \decreases M^{-\alpha}, \quad \textrm{as} \quad M \rightarrow \infty
\end{equation}

\noindent where $\alpha$ is the order of convergence and $\alpha > 0$ \footnote{The notation $a_n \decreases n^{-\alpha}$ denotes that the sequence $\{a_1, a_2, a_3, ..., a_n\}$ decreases with a rate of $-\alpha$ as $n \rightarrow \infty$, and the lower bound of the sequence is $c_1n^{-\alpha}$ and the upper bound is $C_1n^{-\alpha}$ ($c_1n^{-\alpha} \le a_n \le C_1n^{-\alpha}$), where $c_1$ and $C_1$ are positive constants independent of $n$ ($0 < c_1 \le C_1$)}. The value of $\alpha$ is problem dependent and depends on numerous factors, such as, the chosen numerical model, the material model and the smoothness of the random field. By making $M$ sufficiently large, the discretisation error can be reduced to any tolerance value $\epsilon_b$. 

The second term in Eq. (\ref{eq:MC_error_expanded}) is the \textit{sampling error} and represents the variance of the estimator and decays inversely with the number of samples $N$. To ensure that the sampling error is less than a defined tolerance $\epsilon_s$, it is reasonable to determine the number of samples $N$ using Eq. (\ref{eq:MC_num_samples}).

\begin{equation}\label{eq:MC_num_samples}
    \epsilon_s^2 \approx \frac{\mathbb{V}[Q_M]}{N} \quad \therefore N \approx \frac{\mathbb{V}[Q]}{\epsilon_s^2}
\end{equation}

To reduce the total error to a defined tolerance, the number of degrees of freedom $M$ and the number of samples $N$ must both be increased. This can be prohibitively computationally expensive when the cost to compute each sample to the required level of accuracy is high. The cost $C_M$ to compute a single sample of $Q_M$ is dependent on the computational complexity of the solver. The computational cost will grow as follows

\begin{equation}
    C_M \decreases M^{\gamma}
\end{equation}

\noindent for some $\gamma \geq 1$. The rate at which the computational cost grows ($\gamma$) is dependent on a number of factors, such as, the dimension of the problem and the chosen solver (explicit/implicit).

Standard MC estimators are proven to be robust and accurate for many problems, however the slow convergence rate limits applications to problems where the QoI can be computed cheaply. For problems that require the solution of computationally expensive models it is not possible to achieve reasonable estimations in an acceptable time. Different strategies have been examined to accelerate MC estimators, and all are based on the idea of reducing the sampling error.

\subsection{Multilevel Monte Carlo simulation (MLMC)}

The Multilevel Monte Carlo method (MLMC) was introduced by \citet{Giles2008} in 2008, but the first work on multilevel methods was by \citet{Heinrich1999} for parametric integration. Further details on the origins of MLMC are provided in \citet{Giles2015}. Multilevel methods have been widely applied to engineering and scientific problems (i.e. solving partial differential equations). Examples include the computation of the failure probability of composite structures by \citet{Dodwell2021}, the study of the travel time of particles through random heterogeneous porous media by \citet{CrevillenGarcia2017}, and the study of flood risk by \citet{Clare2022}.

\begin{figure}[H]
    \centering
    \includegraphics[width = 0.8\textwidth]{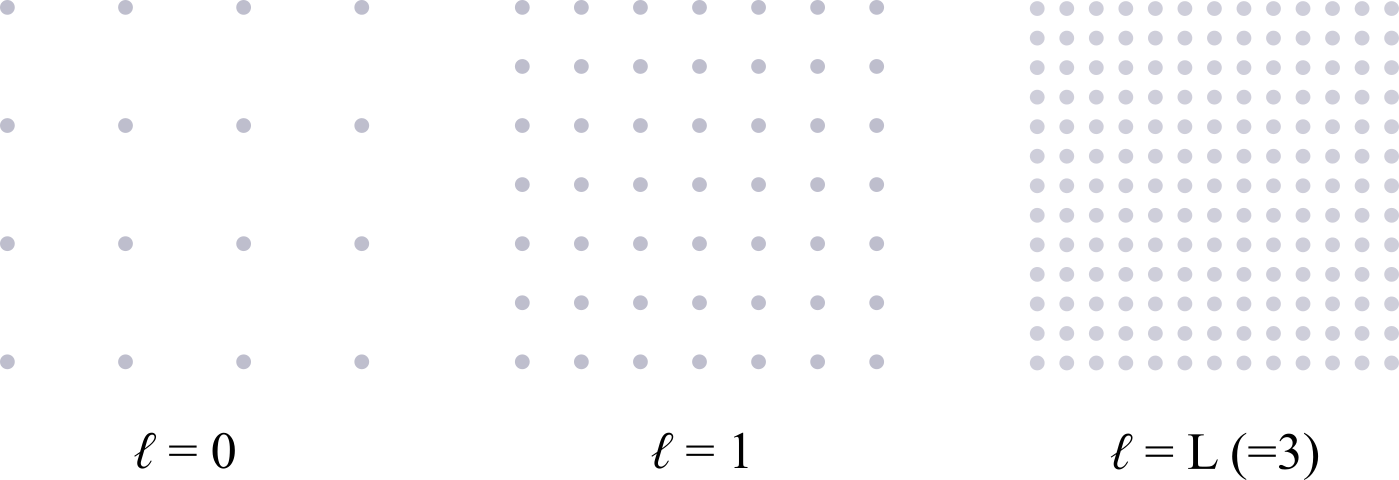}
    \caption{Example of a hierarchy of uniformly refined meshes employed in the MLMC method. Each mesh corresponds to a level $0 \leq \mathbb{\ell} \leq L$ in the multilevel method with $M_0 < \cdots < M_{\ell} < \cdots < M_L$ degrees of freedom. We restrict ourselves to the case of uniform mesh refinement where the node spacing (\(\Delta x\)) is halved every time.}
    \label{fig:nested_meshes}
\end{figure}

The standard MC estimator is too costly as the quantity of interest for every sample must be computed to the level of accuracy required to ensure that the discretisation error is less than a defined tolerance. The key idea of MLMC is to compute a sequence of estimates of the quantity of interest using a hierarchy of nested meshes, as illustrated in Fig. \ref{fig:nested_meshes}. A significant reduction in computational cost can be realised by taking the majority of samples on a coarse mesh (low accuracy but computationally cheap), and taking relatively few samples on a fine mesh (high accuracy but computationally expensive). 

Because of the linearity of the expectation operator, the expected value of $Q$ on the finest mesh ($\mathbb{E}[Q_{M_L}]$) can be expressed as a telescopic sum of the expectation of $Q$ on the coarsest mesh plus a sum of correction terms that account for the difference between evaluations on consecutive mesh levels.

\begin{equation}
    \mathbb{E}[Q_{M_L}] = \mathbb{E}[Q_{M_0}] + \sum_{\ell = 1}^{L}\mathbb{E}[Y_{\ell}]
\end{equation}

\noindent where $Y_{\ell}$ is the discrepancy between the QoI at successive mesh resolutions and is defined as follows

\begin{equation}
    Y_{\ell} = 
    \begin{cases}
        Q_{M_0} & \text{if} \; \ell = 0\\
        Q_{M_{\ell}} - Q_{M_{\ell - 1}} & \text{if} \; 0 < \ell \leq L
    \end{cases}
\end{equation}

The multilevel estimator for $\mathbb{E}[Q]$ is given by Eq. (\ref{eq:MLMC_estimator}).

\begin{equation}\label{eq:MLMC_estimator}
    \hat{Q}^{ML}_{M} = \hat{Q}^{MC}_{M_0, N_0} + \sum^{L}_{\ell=1}\hat{Y}^{MC}_{\ell, N_{\ell}}
\end{equation}

\begin{figure}[H]
    \centering
    \includegraphics[width = 1\textwidth]{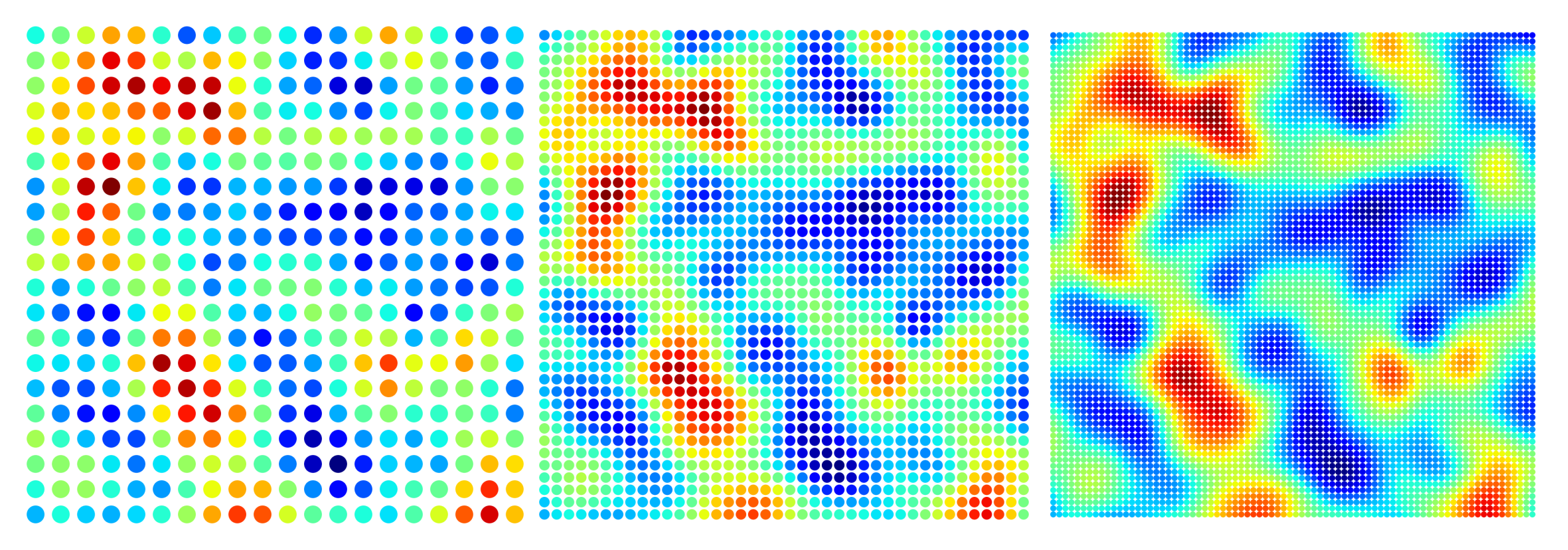}
    \caption{This figure illustrates the same sample $\omega^{(i)}$ of a spatially correlated random field on three mesh levels. Note that the resolution of the coarsest level $\Delta x_0$ must be smaller than the correlation length $l_c$ in the random field. \citet{Cliffe2011} states that the optimal choice for the resolution of the coarsest mesh is such that $\Delta x_0$ is slightly smaller than $l_c$.}
    \label{fig:rf}
\end{figure}

The number of samples $N_{\ell}$ on each level is determined such that the total computational cost of the estimator is minimised for a defined sampling error (see Eq. (\ref{eq:sample_allocation})). It is important to highlight that the same random sample $\omega^{(i)}$ is used to compute the quantity $Q^{(i)}_{M_{\ell}} - Q^{(i)}_{M_{\ell-1}}$, i.e. a coarsened version of the random sample used to compute $Q_{M_{\ell}}^{(i)}$ is used to compute $Q_{M_{\ell-1}}^{(i)}$ (refer to Fig. \ref{fig:rf} for clarification).

As all the expectations $\mathbb{E}[Y_{\ell}]$ are estimated independently, the variance of the multilevel estimator is $\mathbb{V}[\hat{Q}_M^{ML}] = \sum_{\ell=0}^{L}N_{\ell}^{-1}\mathbb{V}[Y_{\ell}]$, where $\mathbb{V}[Y_0] = \mathbb{V}[Q_{M_0}]$. The accuracy of the estimator can be quantified by considering the mean square error.

\begin{equation}\label{eq:MLMC_error}
    e(\hat{Q}_M^{ML})^2 = \underbrace{\left(\mathbb{E}[Q_M - Q]\right)^2}_{bias \; error} + \underbrace{\sum_{l=0}^{L}\frac{\mathbb{V}[Y_{\ell}]}{N_{\ell}}}_{sampling \; error}
\end{equation}

Much like the standard MC estimator, the mean square error is composed of two terms, the bias error and the sampling error. The bias error is exactly the same as in the MC estimator (see Eq. (\ref{eq:MC_error_expanded})), and the number of degrees of freedom on the finest level ($M_L$) must be sufficiently large to satisfy Eq. (\ref{eq:alpha_convergence}), and thus ensuring that the bias error is less than $\epsilon_b$. 

The multilevel estimator is cheaper than the standard MC estimator as the number of samples $N_{\ell}$ on every level can be chosen to ensure that the sampling error is less than $\epsilon_s$, whilst minimising the total computational cost of the estimator. The computational cost of the multilevel Monte Carlo estimator is given by the following

\begin{equation}
    C(\hat{Q}^{ML}_{M}) = \sum_{\ell=0}^{L} N_{\ell} C_{\ell}
\end{equation}

\noindent where $C_{\ell}$ is the cost to compute a single sample of $Y_{\ell}$ on level $\ell \geq 1$ or $Q_{M_0}$ on level 0. Note that taking a sample of $Y_{\ell}$ requires the numerical approximation of $Q$ on two consecutive mesh levels (both $Q_{M_{\ell}}^{(i)}$ and $Q_{M_{\ell-1}}^{(i)}$ must be computed). The determination of the optimal sample allocation is detailed in Section \ref{section:sample_allocation}.

To achieve a RMSE of $\epsilon$, it can be asserted that the multilevel estimator is computationally cheaper than the standard MC estimator due to the significant reduction in variance \cite{Cliffe2011}. As the MLMC estimator is unbiased, the variance of the estimator is equal to

\begin{equation}
    \mathbb{V}[\hat{Q}_M^{ML}] = \sum_{\ell=0}^{L}\frac{\mathbb{V}[Y_{\ell}]}{N_{\ell}}
\end{equation}

The variance of the multilevel estimator is reduced as both numerical approximations $Q_{M_{\ell}}$ and $Q_{M_{\ell-1}}$ converge to $Q$ and consequently

\begin{equation}
    \mathbb{V}[Y_{\ell}] = \mathbb{V}[Q_{M_{\ell}} - Q_{M_{\ell-1}}] \rightarrow 0 \quad \text{as} \quad M_{\ell} \rightarrow \infty
\end{equation}

It is assumed that there exists a $\beta > 0$, where $\beta$ is the order of convergence of the sampling error, such that

\begin{equation}\label{eq:beta_convergence}
    \mathbb{V}[Q_{M_{\ell}} - Q_{M_{\ell-1}}] \decreases M_{\ell}^{-\beta}
\end{equation}

By the central limit theorem, it is clear that fewer samples will be required to accurately approximate the expectation of the difference $Q_{\ell} - Q_{\ell-1}$ as $\ell \rightarrow \infty$. Consequently, the majority of samples will be taken on level 0 (computationally cheap), and relatively few samples will be required at the finest level $L$ (computationally expensive). 

\subsubsection{Error estimation}\label{section:error_estimation}

The aim is to estimate $\mathbb{E}[Q]$ such that the RMSE is below a defined tolerance $\epsilon$, whilst minimising the total computational cost of the estimator $C(\hat{Q}_M^{ML}$). The RMSE, defined by Eq. (\ref{eq:MLMC_error}), is comprised of two parts: (1) the bias error and (2) the sampling error. To ensure that the RMSE is less than $\epsilon$, it is sufficient to bound each term by $\nicefrac{\epsilon^2}{2}$. To estimate the bias error, it is assumed that $M_{\ell}$ is sufficiently large so that the decay in $\big\vert\mathbb{E}[Q_{M_{\ell}} - Q]\big\vert$ is in the asymptotic region and satisfies the following

\begin{equation}\label{eq:alpha_convergence_mlmc}
    \big\vert\mathbb{E}[Q_{M_{\ell}} - Q]\big\vert \decreases M^{-\alpha}
\end{equation}

Following the derivation of \citet{Dodwell2021}, for uniform mesh refinement, where the number of degrees of freedom on level $\ell$ is given by $M_{\ell} \approx m^{\ell}M_0$, the bias error on level $\ell$ can be over-estimated as follows

\begin{equation}\label{eq:bias_error}
    \varepsilon := \big\vert\mathbb{E}[Q_{M_{\ell}} - Q]\big\vert \leq \frac{1}{rm^{\alpha} - 1} \hat{Y}_{\ell, N_{\ell}}^{MC}
\end{equation}

\noindent where $r$ is set to 1. This is equivalent to the assumption that $M_{\ell}$ is sufficiently large so that the decay in $\big\vert\mathbb{E}[Q_{M_{\ell}} - Q]\big\vert$ is in the asymptotic region. The user may wish to select a more conservative values for $r$, for example 0.7 or 0.9. If the bias error is greater than the tolerance, then $M_{L}$ must be increased.

To ensure that the sampling error is less than or equal to the sample tolerance $\epsilon_s$, the following constraint is enforced

\begin{equation}\label{eq:sampling_error}
    \sum_{\ell=0}^L\frac{\mathbb{V}[Y_{\ell}]}{N_{\ell}} \leq \epsilon_s^2
\end{equation}

As the number of samples increases, the variance of the sample mean decreases and hence precision increases. The sample variance is estimated in the standard way \cite{Dodwell2021}

\begin{equation}\label{eq:sample_variance}
    s_{\ell}^2 = \left(\frac{1}{N_{\ell}}\sum_{n=1}^{N_{\ell}}(Y_{\ell}^{n})^2\right) - \left(\hat{Y}_{\ell, N_{\ell}}^{MC}\right)^2 \approx V_{\ell}
\end{equation}

\subsubsection{Sample allocation}\label{section:sample_allocation}

The optimal sample allocation (number of samples per level $N_{\ell}$) is determined by solving a constrained optimisation problem that minimises $C(\hat{Q}^{ML}_{M})$ with respect to $N_{\ell}$, subject to the constraint that the sampling error of the multilevel estimator is less than or equal to the defined tolerance $\epsilon_s$. 

\begin{equation}\label{eq:sample_allocation}
    N_{\ell} = \epsilon_s^{-2}\left(\sum_{\ell=0}^{L} \sqrt{V_{\ell}C_{\ell}} \right) \sqrt{\frac{V_{\ell}}{C_{\ell}}}
\end{equation}

The computational cost of the MLMC estimator grows as follows:

\begin{equation}\label{eq:MLMC_computational_cost}
    C(\hat{Q}_M^{ML}) = \epsilon^{-2}\left(\sum_{\ell=0}^L\sqrt{V_{\ell}C_{\ell}}\right)^2 \decreases \epsilon^{-2-max(0, \; (\gamma - \beta)/\alpha)} \quad \text{as} \quad \epsilon \rightarrow 0
\end{equation}

The rate at which the computational cost grows with respect to the number of degrees of freedom $M$ is given by Eq. (\ref{eq:gamma_convergence}), for some $\gamma \geq 1$. 

\begin{equation}\label{eq:gamma_convergence}
    C_{\ell} \decreases M_{\ell}^{\gamma}
\end{equation}

The reader is referred to \citet{Cliffe2011} and \citet{Giles2015} for a full proof of the MLMC computational complexity theorem with bounds on the RMSE.

\subsubsection{MLMC implementation}\label{section:MLMC_implementation}

Pseudo-code for the adaptive MLMC method is outlined in Algorithm \ref{alg:MLMC}. Optimal values of $L$, $M_{\ell}$ and $N_{\ell}$ are computed \emph{on the fly} from the sample averages and the sample variances of $Y_{\ell}$. We set the number of warm-up samples $N^{\star}$ to be 100. As each sample is independent and there are no shared memory requirements, Algorithm \ref{alg:MLMC} can be trivially parallelised across an unlimited number of independent compute nodes.

\begin{algorithm}[H]
    \caption{Multilevel Monte Carlo Algorithm}\label{alg:MLMC}
    \begin{algorithmic}[1]
    \State Set $L = 0$, $N_{\ell} = N^{\star}$, \texttt{converged} == \texttt{false}

    \While{\texttt{converged} == \texttt{false}}

    \State Take $N^{\star}$ warm-up samples on level $L$
    \State Estimate the sample variance $\mathbb{V}[Y_{\ell}]$ on all levels using Eq. (\ref{eq:sample_variance})
    \State Estimate the optimal number of samples $\hat{N}_{\ell}$ on each level using Eq. (\ref{eq:sample_allocation})
    \State Compute $\hat{N}_{\ell} - N^{\star}$ additional samples on each level
    \State Estimate bias error $\hat{\epsilon}_b$ on level $L$ using Eq. (\ref{eq:bias_error})
    
    \If{$\hat{\epsilon}_b < \epsilon_b$}
        \State \texttt{converged} == \texttt{True}
    \Else 
        \State $L = L + 1$
    \EndIf

    \EndWhile

    \end{algorithmic}
\end{algorithm}




\section{Spatially fluctuating parameter fields}\label{section:random_fields}

To account for the inherent randomness in material properties, spatially fluctuating parameter fields (that vary for every realisation of the parameter field) are employed. In this contribution, we model the spatially fluctuating parameters using random fields. A random field may be thought of as a function that takes random values at every point in the domain. A key property of random fields is that the random variables are spatially correlated. For example, the value of a random variable at adjacent points in space differs less than points that are separated by a large distance. The chosen approach to model this correlation structure can substantially influence the predictive accuracy of the model. The length scale, the maximum distance at which two points are correlated, is another important but ambiguous parameter. 

The efficient and accurate generation of spatially correlated random fields is an important element of the MLMC framework. There are multiple methods for generating spatial random fields and every method has advantages and limitations. Common methods include: turning bands method, spectral method, matrix decomposition method, Karhunen-Lo\`eve expansion (KL expansion) and the moving average method. For a review of the generation of Gaussian random fields and their applications in scientific and engineering problems, the reader is referred to \citet{Liu2019}. For a detailed treatment of the theory of random fields and further applications, the reader is referred to \citet{Hristopulos2020}.

In this contribution, matrix decomposition and KL expansion are employed because of their practical simplicity. The matrix decomposition method generates accurate spatial random fields, but the computational expense is prohibitive for large-scale problems. KL expansion produces less accurate random fields but due to the lower computational cost and ease of implementation, this method was employed for all considered problems. Examining large problems, such as three-dimensional models, is prohibitively expensive. To overcome this issue, the spatial domain can be split into several smaller sub-domains, and a sample of the random field is generated for each sub-domain. A sample for the entire domain is obtained by using an overlapping technique \cite{Panunzio2018, CarvalhoPaludo2019}.

\subsection{Covariance function and length scale}\label{section:covariance_function_and_length_scale}

The value of a random variable at two adjacent points in space is correlated. Conversely, there is negligible correlation between the value of a random variable at two distant points. Many choices exist for the covariance function. Popular choices include Mat\'ern covariance, exponential covariance, Gaussian covariance, spherical covariance and many others \cite{Liu2019}. In this work we have employed an exponential covariance function, as defined by Eq. (\ref{eq:exponential_correlation_function}), where $\rho$ is the correlation coefficient between a random variable at point $\textbf{x}_i$ and $\textbf{x}_j$, $\sigma^2$ is the variance (set to 1), $l_c$ is the correlation length and $\|\textbf{x}_i - \textbf{x}_j\|_2$ is the Euclidean distance between two material points. This form has been selected due to its popularity in the literature. 

\begin{equation}\label{eq:exponential_correlation_function}
    \rho_{ij} = \sigma^2\exp\left(-\frac{\|\textbf{x}_i - \textbf{x}_j\|_2}{l_c}\right)
\end{equation}

There is limited guidance in the literature for selecting a suitable covariance function for different material types. The Joint Commission of Structural Safety (JCSS) Probabilistic Model Code \cite{Vrouwenvelder2001} provides guidance on the probabilistic assessment of concrete structures. The correlation coefficient $\rho_{ij}$ between a random variable at point $\textbf{x}_i$ and $\textbf{x}_j$ is defined by Eq. (\ref{eq:JCSS_correlation_function}).

\begin{equation}\label{eq:JCSS_correlation_function}
    \rho_{ij} = \rho_{ij} + (1 - \rho_{ij})\;\exp\left(-\frac{\|\textbf{x}_i - \textbf{x}_j\|_2}{l_c}\right)
\end{equation}

The covariance function defined in the JCSS probabilistic model code is unusual as the function contains a threshold value for $\rho$. According to the JCSS probabilistic model code the default threshold value is 0.5. To the best of our knowledge, this approach is not seen elsewhere in the literature. By setting the threshold value to 0, the exponential covariance function is obtained. 

\begin{figure}[H]
    \centering
    \includegraphics[width = 1\textwidth]{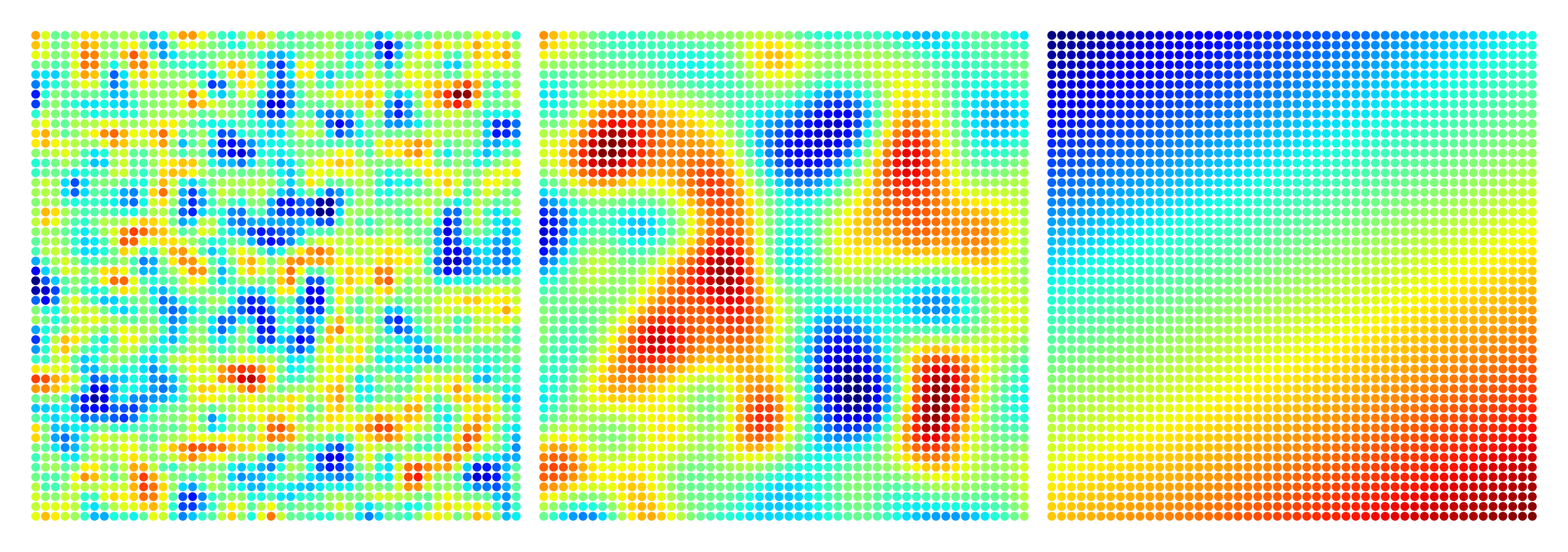}
    \caption{Illustrating the influence of the correlation length $l_c$ on the generated random field (left: $l_c$ = 10 mm, centre: $l_c$ = 30 mm, right: $l_c$ = 5000 mm). The smoothness of the random field increases as $l_c$ increases. Note that the correlation length $l_c$ must be greater than the resolution of the coarsest level $\Delta x_0$ or the value of a random variable at two adjacent points in space will be uncorrelated (white noise).}
    \label{fig:lc}
\end{figure}

The distance over which a correlation exists is determined by the length scale. The correlation length \(l_c\) is a highly uncertain parameter that has a significant influence on the final results. For quasi-brittle materials, \citet{Grassl2009} suggested that the correlation length must, at a minimum, be as large as the fracture process zone (FPZ). For concrete, the size of the FPZ is approximately two to three times the maximum aggregate size \cite{Bazant1989}. The JCSS probabilistic model code recommends a correlation length of 5 m. This value is significantly higher than that found elsewhere in the literature and there is no clear rationale. Fig. \ref{fig:lc} illustrates the influence of the correlation length $l_c$ on the generated random field. It is not the purpose of this paper to examine the correlation length in detail but our studies suggest that a shorter correlation length improves convergence. Note that $l_c$ must be greater than the resolution of the coarsest level $\Delta x_0$ or the value of a random variable at two adjacent points in space will be uncorrelated (white noise). For a detailed examination of the influence of the correlation length $l_c$, the reader is referred to the work of \citet{SyrokaKorol2015}. 

\subsection{Material strength distribution (probability distribution function)}

The chosen material strength distribution plays a key role in the predicted results and convergence of the model. In the literature, normal (Gaussian), log-normal, Gauss-Weibull and Weibull distributions have all been employed for modelling quasi-brittle materials.

The choice of a normal distribution is generally made for convenience as opposed to physical reasons \cite{Bocchini2008}. In particular, material parameters are usually bounded (values must be positive), but negative values are possible when using a normal distribution. Our preliminary studies determined that a normal distribution is not suitable for modelling quasi-brittle materials. This is discussed further in Section \ref{section:problem_2}. The JCSS Probabilistic Model Code recommends that the properties of quasi-brittle materials are modelled using a log-normal distribution \cite{Vrouwenvelder2001}. \citet{VanderHave2015} provides a detailed study of random field generation and the differences between the use of normal and log-normal distributions are explored. The log-normal distribution guarantees that the material parameters are positive, is easy to implement and is widely employed throughout the literature. However, it has been demonstrated that on the scale of a representative volume element (RVE), the probability distribution of strength of quasi-brittle materials is best approximated by a Gaussian distribution onto which a far-left Weibull tail is grafted \cite{Bazant2007b, Le2011, Le2011a}. \citet{Elias2015} and \citet{Elias2020} modelled the size effect in quasi-brittle materials using a lattice discrete particle model (LDPM) and the cumulative distribution function of the random field was assumed to be Gaussian with a left Weibullian tail. The far-left tail of the strength distribution has a huge influence on the the failure load when considering small failure probabilities. For example, for a failure probability of 10$^{-6}$ (structures are generally designed for a failure probability lower than 10$^{-6}$ per lifetime \cite{Melchers2018}), the difference between the failure load and the mean strength will almost double when the strength distribution changes from Gaussian to Weibull (with the same mean and coefficient of variation) \cite{Bazant2007b}. It should be noted that the modelling of quasi-brittle materials is complicated by the transition of the strength distribution from Gaussian to Weibullian as the structure size increases \cite{Le2011}.

It is not the purpose of this paper to examine different strength distributions in detail but we explored the use of a normal, log-normal and Weibull distribution. The Weibull distribution provided the best agreement with experimental data and improved the rate of convergence of the discretisation error. This is discussed further in Section \ref{section:discussion_probability_distribution}.

To easily generate a random field, where the probability distribution function of a material parameter at a given location is a univariate Weibull distribution, we follow the approach of \citet{Rappel2019} and \citet{Rappel2022}. In a Gaussian random field, the probability density function of a material parameter at a given location is a univariate Gaussian distribution. Using the copula theorem and Gaussian fields, different types of univariate marginal distributions can be produced but with the same correlation structure as Gaussian fields. Keeping the Gaussian correlation structure is advantageous as it allows us to draw samples from a Gaussian field and transform the samples into a random field with the desired distribution. The compressive strength $f_c$ is represented as a random field and all other properties are calculated using deterministic relations to $f_c$. Empirical formulas (derived from experimental data) published in \textit{fib} Model Code 2010 \cite{fib2010} are used to determine the Young's modulus $E$, tensile strength $f_t$ and fracture energy $G_F$.


\section{Case studies}\label{section:case_studies}

Two case studies are presented to demonstrate the power of employing the MLMC framework in combination with a peridynamic model. The case studies are selected as examples where uncertainty must be considered to gain a deeper understanding of the physical behaviour. We focus on quasi-brittle materials because the range of experimental data is greater than that of any other material and quasi-brittle materials exhibit a significant size effect. A stochastic model is required for a complete examination of the mechanisms that govern the structural size effect. The first case study examines the structural size effect in a notched beam (Type 2 problem), and the second case study examines the structural size effect in an unnotched beam (Type 1 problem).

The following subsection briefly discuss the structural size effect. For a detailed review of the structural size effect, the reader is referred to \citet{Bazant1998} and \citet{Bazant2000}. 

\subsection{Structural size effect}

Based on the strength-of-materials theory, structural failure is assumed to occur when the maximum stress in a structure exceeds some limiting value of stress that can be determined from small scale tests of representative material samples. Simple fundamental tests such as uniaxial tension, uniaxial compression and flexural tests are used to establish the limiting stress for different loading conditions. This simplistic view does not suffice for quasi-brittle materials \cite{Mier1996}. Quasi-brittle materials exhibit a size effect where large elements fail at lower stresses than geometrically identical but smaller-scale elements.

In brittle and quasi-brittle materials, the size effect can primarily be explained by two mechanisms \cite{Bazant1998,Bazant2000}: (1) release of stored energy (deterministic size effect), and (2) statistical variability in material properties (statistical size effect). The release of stored energy is the principal factor influencing the size effect on structural strength, and the statistical size effect is generally of secondary importance. The deterministic size effect is governed by the size of the fracture process zone (FPZ) (zone of energy dissipation) relative to the size of the structure, and the statistical size effect is a result of the randomness of material properties and defects. The probability that a specimen contains a defect from which failure will initiate increases as the size of the specimen increases. As there is minimal experimental data on the response of very large structures and the safety implications are often much greater, correctly predicting the influence of the statistical size effect is of utmost importance. 

Two types of size effect law are defined: Type 1 - structures with no notches or pre-existing crack (fracture initiates from a smooth surface); Type 2 - structures with a notch. The influence of the size effect on the mean strength of Type 1 and 2 structures is markedly different \cite{Bazant2019}. For Type 1 structures, the size effect has a significant statistical component, for Type 2 structures, the statistical component is minimal.

\citet{Hobbs2022a} previously examined the size effect in quasi-brittle materials using a deterministic bond-based peridynamic model. The model did not consider the spatial variability in material properties and the magnitude of the statistical size effect remains to be established. Due to the high computational expense of peridynamic simulations, examining the statistical size effect was impracticable but the presented framework allows us to overcome the aforementioned issues. \citet{Hobbs2022a} validated the deterministic model against the full set of experimental results published by \citet{Gregoire2013}. This work only considers two members from the test series as the aim of this study is to demonstrate the possible computational savings that can be realised using the MLMC framework, and to demonstrate the importance of examining uncertainty. Future work will use the MLMC framework to examine the full series of tests and provide a comprehensive study of the statistical size effect. 

\subsection{Case study 1: Statistical size effect in quasi-brittle materials (Type 2)}

\begin{figure}[H]
    \centering
    \includegraphics[width = 0.8\textwidth]{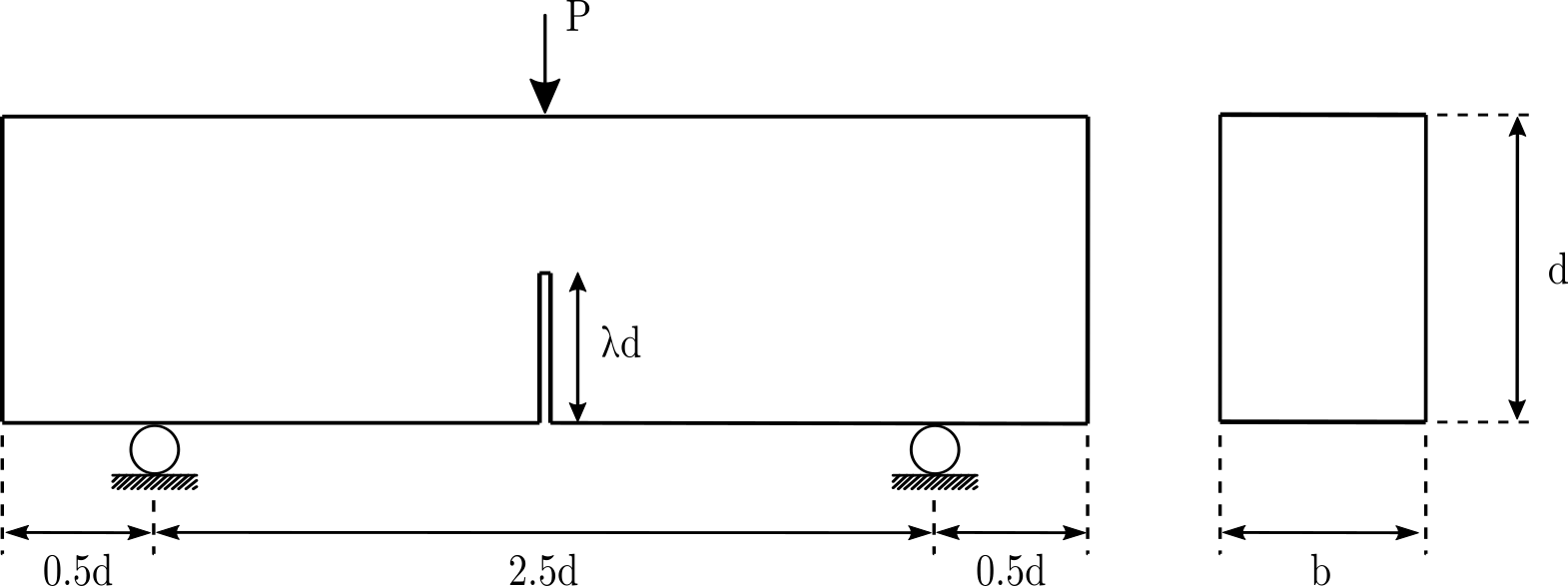}
    \caption{Schematic of the experimental setup (adapted from \citet{Gregoire2013}).}
    \label{fig:mode_I_experimental_setup}
\end{figure}

The first problem that we consider is a notched concrete beam in three-point bending, tested experimentally by \citet{Gregoire2013}. A schematic diagram of the experimental setup is illustrated in Fig. \ref{fig:mode_I_experimental_setup}. The chosen beam (Specimen 3) has the following dimensions: length $l$ = 350 mm; depth $d$ = 100 mm; and thickness $b$ = 50 mm. The span of the member is 250 mm and the depth of the notch is $\lambda$ = 0.5 (half-notched).  

The mean compressive strength $f_{cm,cyl}$ = 42.3 MPa is used to generate a realisation of the random field. The Young's modulus $E$, tensile strength $f_t$ and fracture energy $G_F$ are then estimated using empirical equations \cite{fib2010}. The density of the concrete mixture was 2346 kg/m$^3$ and the maximum aggregate diameter was 10 mm. The correlation length $l_c$ is set to 20 mm. Please refer back to Section \ref{section:covariance_function_and_length_scale} for a discussion of the length scale. The Weibull modulus $m$ is set to 3. This is an uncertain parameter with high sensitivity and a wide range of values can be found in the literature. According to the Weibull theory, the modulus $m$ is a material property that is independent of the geometry and scale of the structure, however \citet{SyrokaKorol2015} found that the Weibull modulus $m$ does depend on the size of the structure and length scale $l_c$.

All the presented results have been obtained using a constant peridynamic horizon $\delta = 3.14\Delta x$ and regular grid spacing. Table \ref{table:discretisation} details the mesh resolution ($\Delta x$) and number of nodes ($M$) for every mesh level ($\ell$). The non-linear softening model (described in Appendix A) was calibrated to fit the experimental results for the smallest unnotched specimen ($d$ = 50 mm). $k$ (the rate of decay) is set to 25 and $\alpha$ (the position of the transition from exponential to linear decay) is set to 0.25. These parameters ($k$ and $\alpha$) are fixed for all test cases.

\begin{table}[H]
    \caption{Mesh level $\ell$, corresponding mesh resolution $\Delta x$ and corresponding number of nodes $M$. The number of nodes (degrees of freedom) on level $\ell$ is given by $M_{\ell}=m^{\ell}M_0$, where $m$ = 4.} 
    \renewcommand{\arraystretch}{1.5} 
    \small
    \centering 

    \begin{tabular}{@{} l *{5}{r} @{}} 
        \toprule
        Level $\ell$ & 0 & 1 & 2 & 3 & 4 \\
        \cmidrule(l){2-6}
        $\Delta x$ (mm) & 10 & 5 & 2.5 & 1.25 & 0.625 \\
        no. nodes $M$ & 350 & 1,400 & 5,600 & 22,400 & 89,600 \\
        \bottomrule
    \end{tabular}
    \label{table:discretisation}
\end{table}

\subsubsection{Results}

We start by taking 100 samples on all levels and estimate $\alpha$, $\beta$ and $\gamma$. The first step is to estimate how the computational cost scales as $M_{\ell}$ increases. The time to compute each sample is recorded and it is determined that the computational costs grows linearly. The computational cost is given by Eq. (\ref{eq:type_2_gamma_convergence}), where $\gamma$ = 1. Note that the performance of \texttt{PeriPy} scales linearly with the number of nodes $\therefore \gamma = 1$ \cite{Boys2021a}. 

\begin{equation}\label{eq:type_2_gamma_convergence}
    C_{\ell} \decreases M_{\ell}^{1}
\end{equation}

\begin{figure}[H]
    \centering
    \includegraphics[width = 1\textwidth]{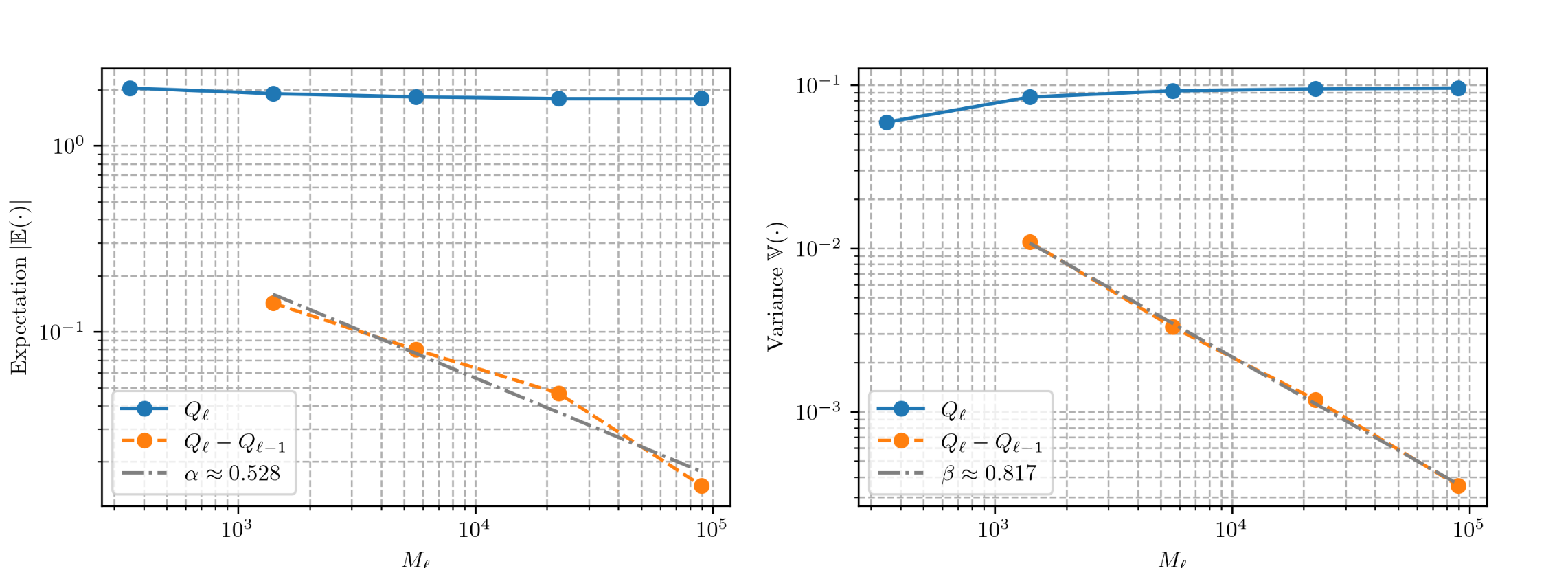}
    \caption{Expectation (left) and variance (right) of $Q_{\ell}$ and $Y_{\ell} = Q_{\ell} - Q_{\ell - 1}$ against degrees of freedom $M_{\ell}$ for problem 1. Taking 100 samples on every level, $\alpha$ is estimated to be 0.528 and $\beta$ is estimated to be 0.817.}
    \label{fig:B3_HN_convergence}
\end{figure}

The next step is to estimate the parameters $\alpha$ and $\beta$ for the QoI, which is taken to be the peak load. Fig. \ref{fig:B3_HN_convergence} illustrates the log-log plots of the estimated means and variances of $Q_{\ell}$ and $Y_{\ell} = Q_{\ell} -  Q_{\ell - 1}$, for $\ell = 0, ..., 4$, with respect to the number of degrees of freedom $M_{\ell}$ on each level. The rate of convergence of the discretisation error is given by Eq. (\ref{eq:type_2_alpha_convergence}), where $\alpha$ is approximately 0.528.

\begin{equation}\label{eq:type_2_alpha_convergence}
    \mathbb{E}[Y_{\ell}] \decreases M_{\ell}^{-0.528}
\end{equation}

The rate of convergence of the sampling error is given by Eq. (\ref{eq:type_2_beta_convergence}), where $\beta$ is approximately 0.817.

\begin{equation}\label{eq:type_2_beta_convergence}
    \mathbb{V}[Y_{\ell}] \decreases M_{\ell}^{-0.817}
\end{equation}

\begin{figure}[H]
    \centering
    \includegraphics[width = .6\textwidth]{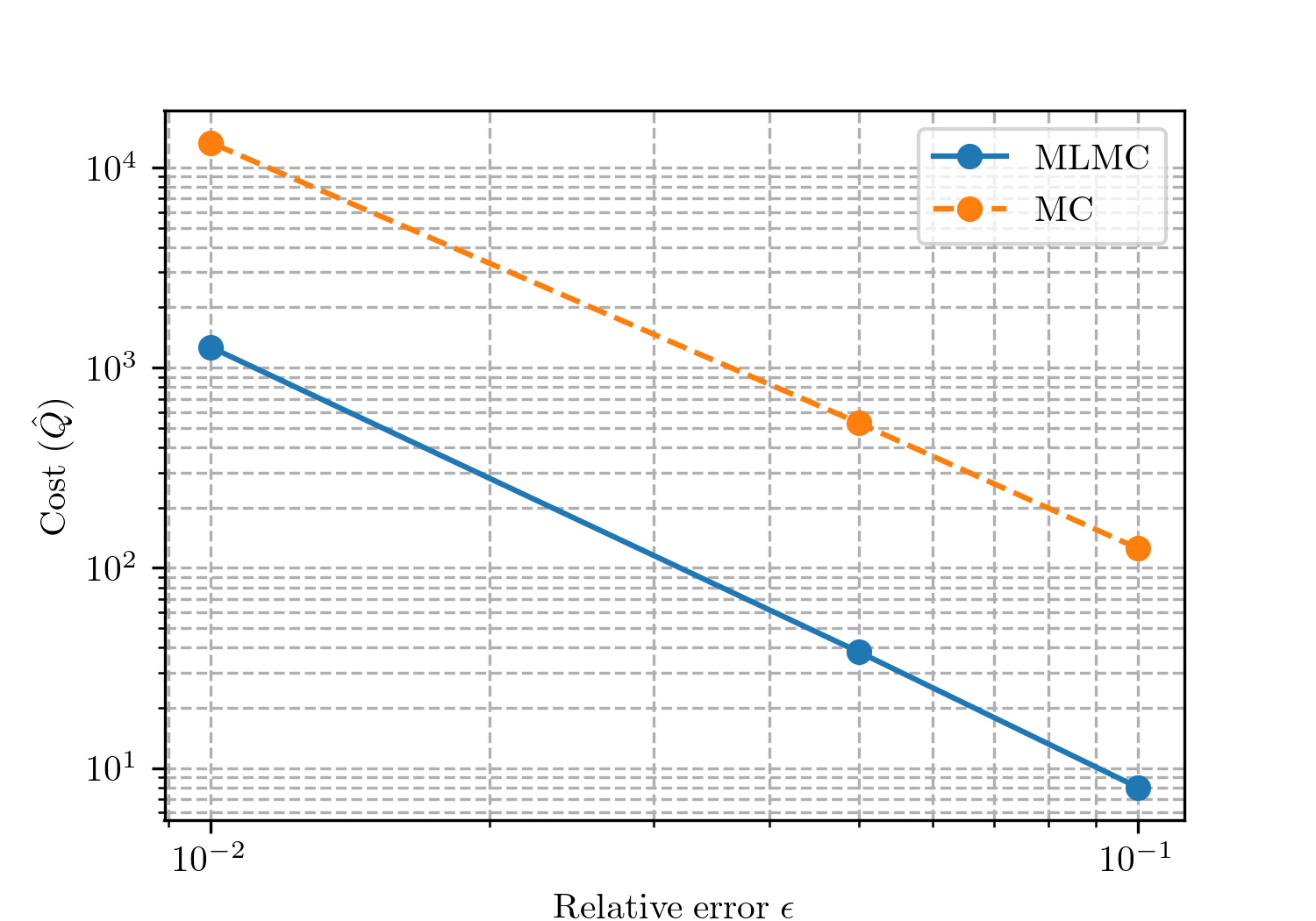}
    \caption{Computational cost (in minutes) against $\epsilon$ for the multilevel estimator (cost $\decreases \; \epsilon^{-2.35}$) and the standard MC estimator (cost $\decreases \; \epsilon^{-3.89}$)}
    \label{fig:Cost_vs_epsilon_Type_2}
\end{figure}

Using the estimated values of $\alpha$, $\beta$ and $\gamma$, Eq. (\ref{eq:MLMC_computational_cost}) predicts that the cost of the MLMC simulations will grow proportionally to $\epsilon^{-2.35}$, whilst the cost of the standard MC simulations will grow proportionally to $\epsilon^{-3.89}$. Fig. \ref{fig:Cost_vs_epsilon_Type_2} illustrates the computational cost (in minutes) of the multilevel estimator and the standard MC estimator against relative error $\epsilon$. Table \ref{table:type_2_sample_allocation} presents the optimal number of samples $N_{\ell}$ across the 5 mesh levels for different values of sampling tolerance ($\epsilon_s$ = 10, 50 and 100 N), plus the number of samples required when using the standard MC estimator ($N$). Using Eq. (\ref{eq:bias_error}), the bias error is found to be approximately 0.75 N. 

\begin{table}[H]

    \caption{Specimen 3 ($\lambda$ = 0.5) - Sample allocation for different sampling tolerances $\epsilon_s$. The sampling tolerance $\epsilon_s$ is given in Newtons (N). $N_{\ell}$ is computed using Eq. (\ref{eq:sample_allocation}) and $N$ is computed using Eq. (\ref{eq:MC_num_samples}). The bias error is approximately 0.75 N}
    \centering 
    
    \begin{threeparttable}
    
    \renewcommand{\arraystretch}{1.5} 
    \footnotesize

        \begin{tabular}{@{} l *{6}{c} @{}} 
            \toprule
            \multirow{2}{*}{$\epsilon_s$ (N)} & \multicolumn{5}{c}{no. samples $N_{\ell}$} & \multirow{2}{*}{$N$\tnote{1}} \\
            \cmidrule(lr){2-6}
            \addlinespace[-0.6ex]
            & 0 & 1 & 2 & 3 & 4 & \\
            \midrule
            100 & 21 & 5 & 2 & 0 & 0 & 9 \\
            50 & 86 & 22 & 8 & 2 & 0 & 38 \\
            10 & 2,152 & 567 & 204 & 67 & 18 & 955 \\
            \bottomrule
        \end{tabular}

    \begin{tablenotes}
        \item[1] Number of samples required when using the standard Monte Carlo estimator. Note that all samples are computed on level 4 ($\ell = 4$). 
    \end{tablenotes}

    \end{threeparttable}
    \label{table:type_2_sample_allocation}
\end{table}

\begin{table}[H]

    \caption{Specimen 3 ($\lambda$ = 0.5) - Computational cost (in minutes) for the multilevel estimator and the standard Monte Carlo estimator.} 
    \centering 
        
    \renewcommand{\arraystretch}{1.5} 
    \footnotesize

    \begin{tabular}{@{} l *{3}{c} @{}} 
        \toprule
        \multirow{2}{*}{$\epsilon_s$ (N)} & \multicolumn{2}{c}{Cost (min)} & \multirow{2}{*}{Speed-up} \\
        \cmidrule(lr){2-3}
        \addlinespace[-0.6ex]
        & MLMC & MC & \\
        \midrule
        100 & 8 & 125 & 16 \\
        50 & 38 & 529 & 14 \\
        10 & 1,265 & 13,290 & 10.5 \\
        \bottomrule
    \end{tabular}

    \label{table:type_2_computational_cost}
\end{table}

Whilst the main aim of this work has been to demonstrate the computational savings of the MLMC framework, we have purposely selected examples where a deeper understanding of the physical behaviour can be gained by considering uncertainty. For Type 2 problems, the difference between the deterministic strength and the mean stochastic strength is expected to be small, as the influence of the variability in material properties is lessened due to the large stress concentrations that occur at the notch tip. This was observed numerically by \citet{Elias2015} who found that considering spatial variability in material properties does not significantly influence the mean failure load, but does lead to an increase in the variance of the structural response. Using the finest mesh ($\ell$ = 4), the deterministic model of \citet{Hobbs2022a} predicts that the specimen will fail at approximately 1800 N. Setting the sampling tolerance $\epsilon_s$ to 10 N, the mean stochastic strength is predicted to be approximately 1790 N. Note that the bias error is approximately 0.75 N. The observed results are in agreement with theory, which predicts that the difference between the deterministic strength and mean stochastic strength will be small \cite{Bazant1984,Bazant2019}. Note that the experimental failure load ranged between 1580 N and 1710 N. 

\subsection{Case study 2: Statistical size effect in quasi-brittle materials (Type 1)}\label{section:problem_2}

The second problem that we consider is an unnotched concrete beam in three-point bending, tested experimentally by \citet{Gregoire2013}. We consider Specimen 3 (illustrated in Fig. \ref{fig:mode_I_experimental_setup}) again but with no notch ($\lambda$ = 0). Beyond demonstrating the computational savings that can be achieved using the MLMC framework, the presented example provides insight into the following areas:

\medskip

\noindent \textbf{Statistical size effect} - \citet{Hobbs2022a} showed that a deterministic bond-based model accurately captures the structural size effect for Type 2 (notched) problems, but fails to capture the correct response for Type 1 (unnotched) problems. This was expected as it is well known that the randomness of material properties has a significant effect on the structural strength of Type 1 problems \cite{SyrokaKorol2013, Elias2020}. In Type 1 problems, the volume of highly stressed material is much larger than that observed in Type 2 problems, and the probability that a defect is present in the stressed region is consequently higher. In Type 2 problems, the presence of a notch results in a localised region of highly stressed material, and the influence of randomness in material properties is consequently lessened. It is expected that the inclusion of statistical variability in the material properties will improve the predictive accuracy of the peridynamic model.

\medskip

\noindent \textbf{Convergence} - \citet{Hobbs2021a} demonstrated that a deterministic bond-based model fails to converge for Type 1 problems (the predicted strength is coupled with the mesh resolution). It was hypothesised that accounting for randomness in the material properties is required to initiate the localisation of damage and improve convergence. 

\subsubsection{Results}

\begin{figure}[H]
    \centering
    \includegraphics[width = 1\textwidth]{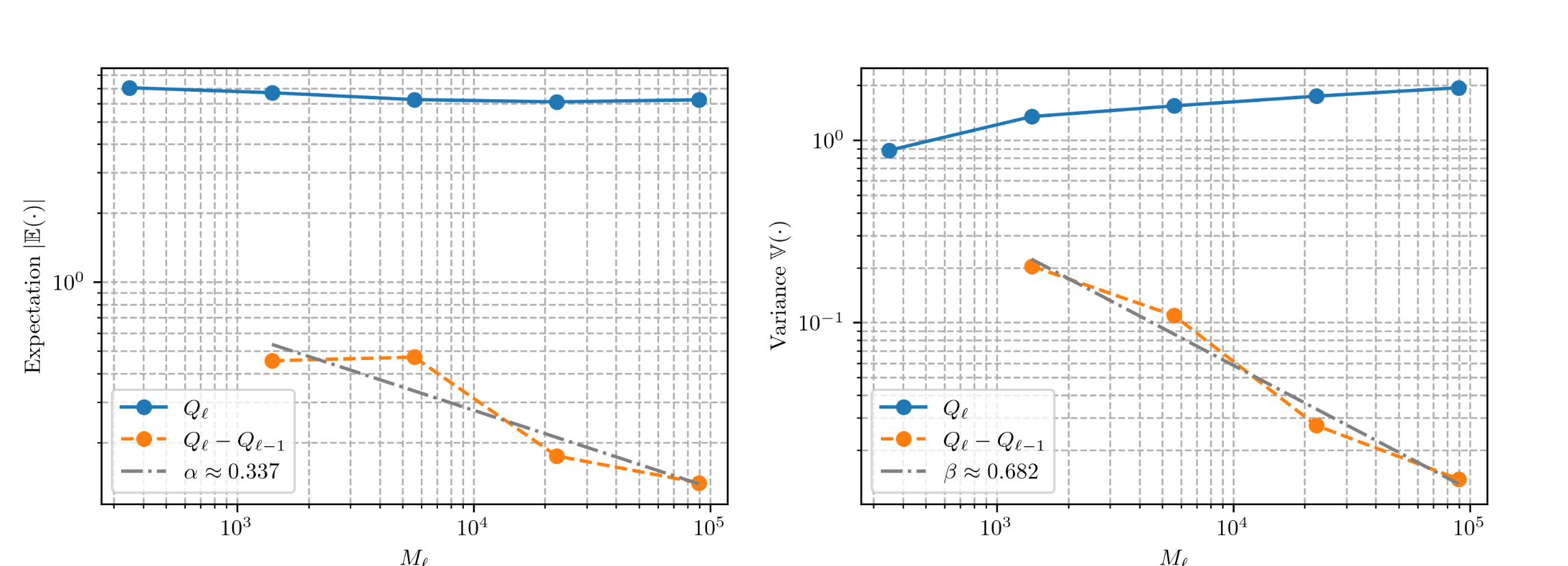}
    \caption{Expectation (left) and variance (right) of $Q_{\ell}$ and $Y_{\ell} = Q_{\ell} - Q_{\ell - 1}$ against degrees of freedom $M_{\ell}$ for problem 2. Taking 100 samples on every level, $\alpha$ is estimated to be 0.337 and $\beta$ is estimated to be 0.682.}
    \label{fig:B3_UN_convergence}
\end{figure}

Again we start by taking 100 samples on all levels and estimate $\alpha$, $\beta$ and $\gamma$. As per the previous example, the computational cost grows linearly ($\gamma$ = 1). Taking 100 samples on every level, $\alpha$ is estimated to be 0.337 and $\beta$ is estimated to be 0.682 (refer to Fig. \ref{fig:B3_UN_convergence}). The rate of convergence of the discretisation error and sampling error is slower than that observed in problem 1 (Type 2).

Using the estimated values of $\alpha$, $\beta$ and $\gamma$, Eq. (\ref{eq:MLMC_computational_cost}) predicts that the cost of the MLMC simulations will grow proportionally to $\epsilon^{-2.94}$, whilst the cost of the standard MC simulations will grow proportionally to $\epsilon^{-4.97}$. Table \ref{table:type_1_sample_allocation} presents the optimal number of samples $N_{\ell}$ across the mesh levels for different values of sampling tolerance ($\epsilon_s$ = 10, 50 and 100 N), plus the number of samples required when using the standard MC estimator ($N$). Due to the higher variance of the estimator, the number of samples required is considerably higher than that required for the Type 2 problem. Type 1 problems are subject to a high degree of natural variability and consequently the computational cost is higher as significantly more samples are required. 

\begin{table}[H]

    \caption{Specimen 3 ($\lambda$ = 0) - Sample allocation for different sampling tolerances $\epsilon_s$. The sampling tolerance $\epsilon_s$ is given in Newtons (N). $N_{\ell}$ is computed using Eq. (\ref{eq:sample_allocation}) and $N$ is computed using Eq. (\ref{eq:MC_num_samples}). Due to the higher variance of the estimator, the number of samples required is considerably higher than that required for the Type 2 problem (refer to Table \ref{table:type_2_sample_allocation}).} 
    \centering 
    
    \begin{threeparttable}
    
    \renewcommand{\arraystretch}{1.5} 
    \footnotesize

        \begin{tabular}{@{} l *{6}{c} @{}} 
            \toprule
            \multirow{2}{*}{$\epsilon_s$ (N)} & \multicolumn{5}{c}{no. samples $N_{\ell}$} & \multirow{2}{*}{$N$\tnote{1}} \\
            \cmidrule(lr){2-6}
            \addlinespace[-0.6ex]
            & 0 & 1 & 2 & 3 & 4 & \\
            \midrule
            100 & 406 & 119 & 57 & 15 & 5 & 193 \\
            50 & 1,624 & 479 & 230 & 63 & 22 & 775 \\
            10 & 40,624 & 11,982 & 5,765 & 1,586 & 551 & 19,388 \\
            \bottomrule
        \end{tabular}

    \begin{tablenotes}
        \item[1] Number of samples required when using the standard Monte Carlo estimator. Note that all samples are computed on level 4 ($\ell = 4$). 
    \end{tablenotes}

    \end{threeparttable}
    \label{table:type_1_sample_allocation}
\end{table}

\begin{table}[H]

    \caption{Specimen 3 ($\lambda$ = 0) - Computational cost (in minutes) for the multilevel estimator and the standard Monte Carlo estimator. The computational cost of the multilevel estimator grows proportionally to $\epsilon^{-2.94}$ and the cost of standard Monte Carlo estimator grows proportionally to $\epsilon^{-4.97}$.} 
    \centering 
        
    \renewcommand{\arraystretch}{1.5} 
    \footnotesize

    \begin{tabular}{@{} l *{3}{c} @{}} 
        \toprule
        \multirow{2}{*}{$\epsilon_s$ (N)} & \multicolumn{2}{c}{Cost (min)} & \multirow{2}{*}{Speed-up} \\
        \cmidrule(lr){2-3}
        \addlinespace[-0.6ex]
        & MLMC & MC & \\
        \midrule
        100 & 290 & 2,686 & 9.3 \\
        50 & 1,197 & 10,785 & 9 \\
        10 & 29,998 & 269,816 & 9 \\
        \bottomrule
    \end{tabular}

    \label{table:type_1_computational_cost}
\end{table}

By including uncertainty in the material properties, the bond-based model converges for Type 1 problems ($\alpha \approx 0.337$). This is the first time that this behaviour has been demonstrated, but the convergence behaviour is significantly worse than that observed for Type 2 problems ($\alpha \approx 0.528$) and the bias error remains large. Using Eq. (\ref{eq:bias_error}), the bias error is estimated to be approximately 200 N.

Initial results using a normal distribution found that the model failed to converge for Type 1 problems. It was established that extreme values in the left-tail of the material strength distribution are required to initiate the localisation of damage and eliminate problems of mesh dependence. If there are no imperfections in the highly stressed region then the damage fails to localise correctly and the results exhibit a strong mesh dependency. Even when using a Weibull distribution, there will be a small number of samples where the damage fails to localise and this has a detrimental impact on the estimated convergence rate $\alpha$. The convergence rate $\alpha$ can be improved by employing a material strength distribution that is skewed towards the left (e.g. Weibull distribution with a low Weibull modulus) but this might not be physically realistic for the considered problem.

As the size of a structure increases, so does the probability that a defect will be present from which a fracture will initiate. \citet{SyrokaKorol2013} determined numerically that the deterministic and mean stochastic strength start to diverge when the beam depth is greater than 50-60 mm. Specimen 3 is 100 mm deep and the magnitude of the statistical size effect is expected to be non-negligible. Setting the sampling tolerance $\epsilon_s$ to 50 N, the mean stochastic strength is estimated to be approximately 6250 N. Note that the bias error is approximately 200 N. Using the finest mesh ($\ell$ = 4), the deterministic model predicts that the specimen will fail at approximately 9200 N. The experimental failure load ranged between 7620 N and 8770 N. The numerical results are consistent with the theory, i.e., the difference between the deterministic strength and mean stochastic strength is much larger than that observed for Type 2 problems. However, the deterministic model does not converge for Type 1 problems and the prediction of strength is therefore unreliable, and a rigorous comparison is not possible. 


\section{Cumulative distribution function}\label{section:cdf}

The objective of the multilevel framework is to estimate the expectation of an output variable, in this case, the peak load. However, for many industrial applications, engineers are more concerned with the cumulative distribution function (CDF) of the output variable. For example, an engineer might be interested in the probability that an output variable exceeds a specific value, or as demonstrated here, the probability that a structure will fail at a load less than or equal to a specific value. Computing the CDF is complicated as the multilevel method provides relatively few values on the finest mesh. \citet{Gregory2017} recently outlined a method that makes it possible to obtain the CDF using samples obtained on multiple mesh levels and we follow the same approach. The reader is also referred to \citet{Clare2022} for further information.

\citet{Gregory2017} use the inverse transform sampling method to determine an approximation of the inverse CDF $F^{-1}(u)$, where $u \sim \mathcal{U}(0, 1)$. If the CDF $F$ of a random variable $X$ is strictly increasing and absolutely continuous, then there exists a unique value $x \equiv F^{-1}(u)$ for which $F(x) = u$. By sorting the samples $\{X^i\}_{i=1,..., N} \sim F$ in ascending order such that $X^1 < X^2 < \cdots < X^N$, a simple consistent estimate for an evaluation to the quantile function of the distribution with CDF $F$ is 

$$\hat{F}^{-1}(u) = X^{[N \times u]}$$

This is a consistent estimate because it converges in probability to $F^{-1}(u)$ as $N \rightarrow \infty$. The inverse CDF for the multilevel approximation is then given by

$$\hat{F}^{-1}_L(u) = R(X)_0^{[N_0 \times u]} + \sum^L_{i=1}\left(R(X)_{\ell}^{[N_{\ell} \times u]} - R(X)_{\ell-1}^{[N_{\ell-1} \times u]}\right)$$

\noindent where $R(X)^i_{\ell}$ is the $i$th ordered statistic of $X_{\ell}$ on level $\ell$. The CDF of the multilevel approximation is then given by 

$$\hat{F}(x) = \frac{1}{N}\sum^{N}_{i=1}\mathbbm{1}_{X_i \leq x}$$

\noindent where $\mathbbm{1}$ is the indicator function. For a detailed description of this method we refer the reader to \citet{Gregory2017}.

\begin{figure}[H]
    \centering
    \includegraphics[width = 1\textwidth]{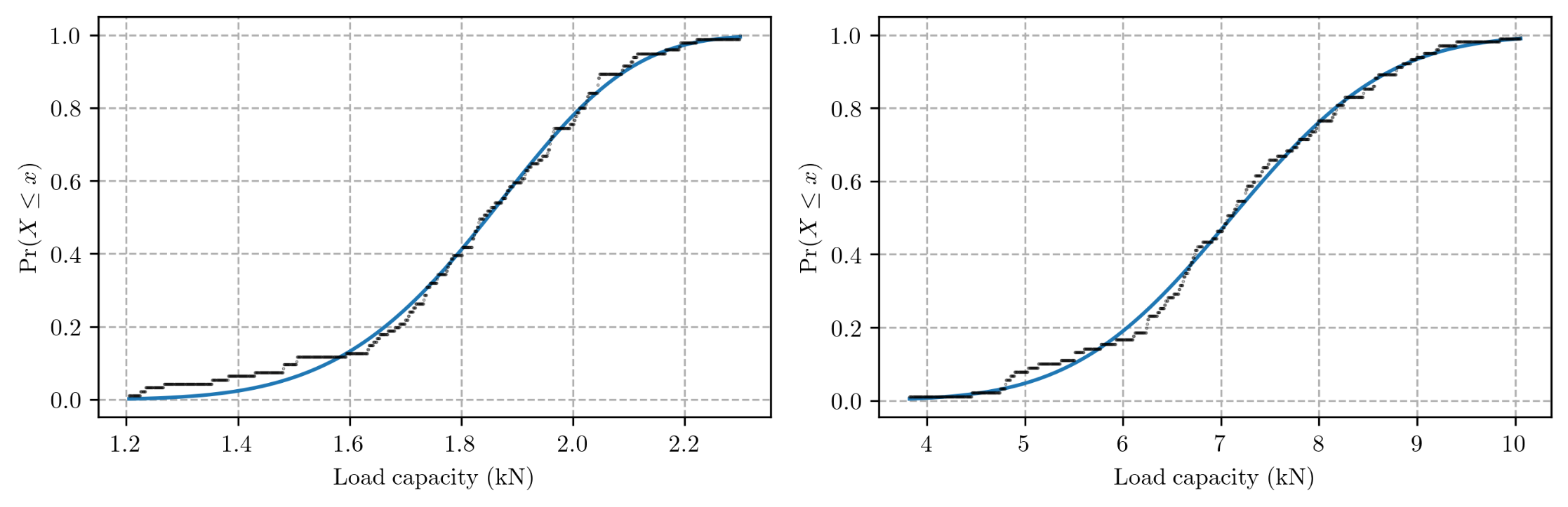}
    \caption{The cumulative distribution function of strength (load capacity) for the Type 2 problem (left) and the Type 1 problem (right). The $y$-axis represents the percentage of the population that will fail at a load less than or equal to $x$.}
    \label{fig:cdf}
\end{figure}

The empirical CDF and fitted Weibull CDF for the Type 1 and Type 2 problem is illustrated in Fig. \ref{fig:cdf}. The $y$-axis represents the percentage of the population that will fail at a load less than or equal to $x$. The purpose of this figure is to demonstrate that is possible to use the multilevel framework to obtain the CDF and a discussion of the results is beyond the scope of this paper. The reader is referred to \citet{Bazant2017} for further information on the computation of CDFs of structural strength and fatigue lifetime. Ideally we would include the CDF obtained using standard MC but this was not possible due to the high computational cost of obtaining samples on the finest level.  

\section{Discussion}\label{section:discussion}

\subsection{Statistical size effect}

A key aim of this study was to select case studies where uncertainty must be considered to gain a comprehensive understanding of the physical behaviour. We focussed our studies on the structural size effect in quasi-brittle materials. \citet{Bazant1991a} stated that the correct modelling of the size effect on material strength should be adopted as the basic criterion of acceptability of any model. The results demonstrate that a bond-based peridynamic model can be used to examine both the statistical and deterministic component of the structural size effect. The intention of this study was never to provide a detailed examination of the statistical size effect, and further studies on a wider range of problems are required to improve confidence in the models predictive capabilities. 

By employing the presented MLMC framework, studying the statistical component of the structural size effect using a peridynamic model becomes computationally feasible. Future work aims to employ the presented MLMC framework to study the full series of tests published by \citet{Gregoire2013} and provide a detailed examination of influential factors, such as the shape of the material strength distribution and the correlation length $l_c$. \citet{Grassl2009} state that the ratio of the correlation length $l_c$ to the size of the fracture process zone (FPZ) is the main parameter that influences the statistical size effect. 

\subsection{Convergence}

Numerical results should be independent of the mesh resolution. This is a basic test of the adequacy of any numerical model. To the best of the authors knowledge, \citet{Hobbs2021a} was the first to consider the effect of mesh refinement ($\delta$-convergence) on the predicted peak load and load-deflection response for Type 1 and 2 problems. Hobbs found that a deterministic bond-based peridynamic model fails to converge for Type 1 problems. Note that \citet{Niazi2021} also published a convergence study that considered the complete structural response. The study of \citet{Niazi2021} is limited as Type 1 problems were not considered. 

The results in this study confirm that, as previously hypothesised, a source of randomness must be introduced to trigger the localisation of damage in Type 1 specimens and eliminate problems of mesh dependence that occur in peridynamic models. \citet{Niazi2021} reported that the convergence behaviour is improved by randomly deleting 1\% of all bonds, as first suggested by \citet{Chen2019b}. Whilst the method of \citet{Chen2019b} is computationally cheap and does improve convergence behaviour, it is an oversimplified approach that lacks a robust theoretical basis (heuristic) and does not consider the spatial correlation of material properties. \citet{Jones2020} note that these methods are generally used to avoid problems related to symmetry, and they do not attempt to capture the true material behaviour by implementing an experimentally measured probability distribution of material properties. 

\subsection{Length scales}

The correlation length $l_c$ was set to be 20 mm for all considered problems. This value was selected after running a number of preliminary simulations. However, the aim of this contribution was not to identify the parameters that describe the spatial fields. It is important to note that a theoretically grounded probabilistic framework based on Bayesian inference (see \cite{Rappel2019, Rappel2022, Rappel2020}) is essential to identify the parameters of the spatial fields (e.g. length scale $l_{c}$) rigorously. 

Furthermore, the interaction between the two length scales (peridynamic horizon $\delta$ and the correlation length $l_c$ in the random field) requires further examination. It remains uncertain how the ratio of the two length scales influences the predictive accuracy of the model. 

\subsection{Probability distribution}\label{section:discussion_probability_distribution}

The material strength distribution plays an important role in the predicted results and convergence of the model. Three distribution were considered (normal, log-normal and Weibull) and it was determined that the Weibull distribution provides the best predictions of mean strength for quasi-brittle materials. This was expected and has been extensively discussed in the literature. A more novel observation is that the selected probability distribution influences the convergence rate of the bias error. Extreme values in the left-tail are required to initiate the localisation of damage and eliminate problems of mesh dependence. Note that the model failed to converge for Type 1 problems when using a normal distribution. 

\subsection{Model calibration}

Many of the model parameters are impossible to determine exactly and are subject to significant uncertainties, for example: the length scale $l_c$ and the Weibull modulus (shape parameter). Future work will examine the integration of the multilevel method with experimental data in a Bayesian setting to quantify modelling uncertainties as proposed by \citet{Dodwell2015, Dodwell2019a}. This will be an important step in the validation of peridynamic models, enabling the identification of model discrepancy and measurement bias, and providing better estimates of model parameters.


\section{Conclusions}\label{section:conclusions}

Peridynamic models are computationally expensive, thus preventing the use of standard Monte Carlo methods for the assessment of uncertainties in model outputs propagated from uncertain inputs. The aim of this study was to demonstrate the possible computational savings that can be realised using the MLMC framework. The results show a speed-up factor of 16$\times$ over a standard Monte Carlo estimator, enabling the forward propagation of uncertain parameters in a computationally expensive peridynamic model. Beyond demonstrating the computational savings that can be achieved using the multilevel framework, the results presented in this paper are of interest for two further reasons:

\begin{enumerate}
    \item Deterministic bond-based models suffer from a strong mesh dependency when simulating Type 1 problems. It has been demonstrated that by including uncertainty in the material properties, the bond-based peridynamic model converges for both Type 1 and Type 2 problems. The need to consider uncertainty is essential for robust and accurate predictions. Furthermore, the multilevel method provides an estimate of the discretisation error, thus improving the interpretability of numerical predictions.  
    \item A secondary aim was to select case studies where uncertainty must be considered to gain a comprehensive understanding of the physical behaviour. We examined the structural size effect in quasi-brittle materials as the random variability of material properties is known to play an important role. \citet{Bazant1991a} stated that the correct modelling of the size effect on material strength should be adopted as the basic criterion of acceptability of any model. The results demonstrate that a bond-based peridynamic model can be used to study the statistical size effect but further studies on a wider range of problems are required to improve confidence in the models predictive capabilities. Future work will consider the full series of tests published by \citet{Gregoire2013} and provide a detailed study of the statistical size effect. 
\end{enumerate}

We have motivated the use of the MLMC framework by studying the statistical size effect in quasi-brittle materials. But forward uncertainty quantification is equally important for cases where a high degree of reliability is required, as is common in many aerospace and power generation applications. 


\section*{Acknowledgements}

TD and MH are funded through a Turing AI Fellowship, United Kingdom (2TAFFP\textbackslash100007). The authors would like to acknowledge the use of the University of Exeter High-Performance Computing (HPC) facility in carrying out this work.


\section*{Data access statement}


\section*{Appendix A. Non-linear softening model (two-dimensional case)}

\begin{figure}[H]
    \centering
    \includegraphics[width = 0.55\textwidth]{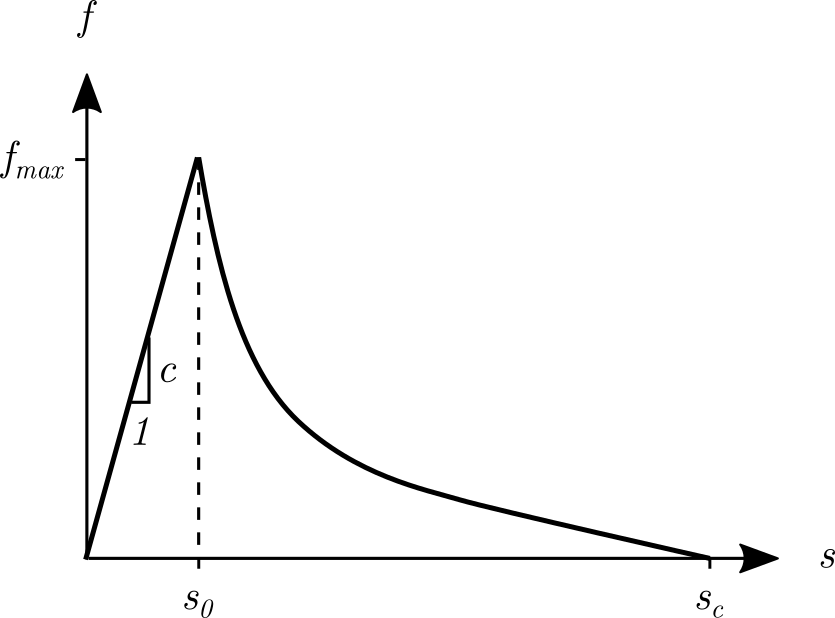}
    \caption{Non-linear damage model ($f$-$s$). The force-stretch relationship is described by an exponentially decaying model with a linear term that forces the curve to intersect with the horizontal axis at $s_c$.}
    \label{fig:Decaying_exponential}
\end{figure}

We employ the non-linear model introduced by \citet{Hobbs2021a} (Section 5.3.4) and employed in \citet{Hobbs2022a}. The model is illustrated in Fig. \ref{fig:Decaying_exponential}. Here we provide the derivation for the two-dimensional plane stress and plane strain case. Note that the Poisson's ratio $\nu$ is limited to 1/3 for plane stress problems, and 1/4 for plane strain and 3D problems. 

\begin{figure}[H]
    \centering
    \includegraphics[width = 0.6\textwidth]{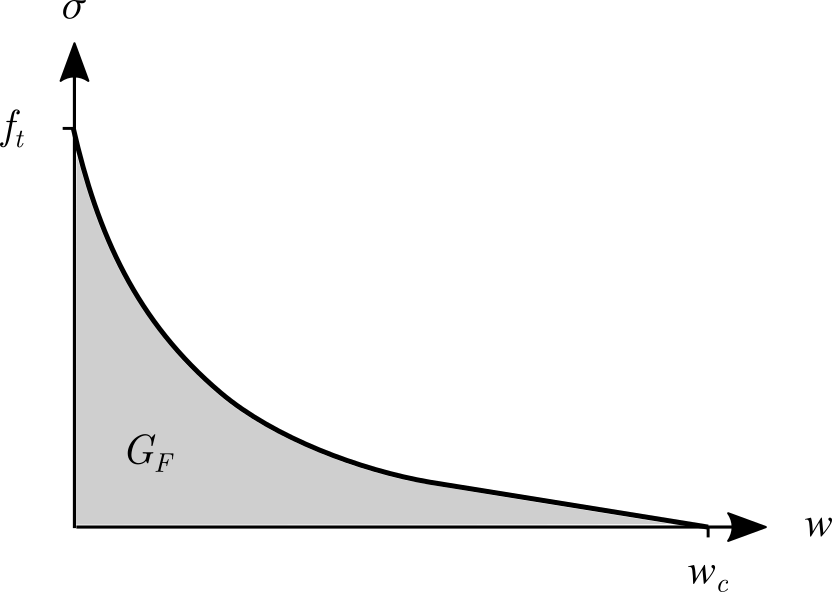}
    \caption{Non-linear tension softening constitutive law ($\sigma$-$w$). Adapted from Fig. 6.1 in \citet{Hordijk1991}.}
    \label{fig:Non-linear_softening_law}
\end{figure}

It is assumed that the behaviour of an individual peridynamic bond is consistent with the stress-crack width ($\sigma$-$w$) relationship illustrated in Fig. \ref{fig:Non-linear_softening_law}. The constitutive relationship was derived for quasi-brittle materials from the experimental work of \citet{Cornelissen1986}. Note that the area under the $\sigma$-$w$ curve is a measure of the material fracture energy $G_F$.

The $\sigma$-$w$ relationship is described by an exponentially decaying model with a term that forces the curve to intersect with the horizontal axis at $w_c$. If the softening relationship is asymptotic with the horizontal axis, and thus never intersects, a unique value for the critical stretch of a bond cannot be determined.

The scalar bond force $f$ is defined by Eq. (\ref{eq:exp_scalar_force}), where $c_d$ is the stiffness of a damaged bond, $s$ is the bond stretch, $c$ is the stiffness of an undamaged bond and $d$ is a softening parameter.

\begin{equation}\label{eq:exp_scalar_force}
    f = c_ds \quad c_d = c(1 - d)
\end{equation}

The bond stiffness constant for different problem types is defined by Eq. (\ref{eq:bond_stiffness_constant}), where $t$ is the thickness of the domain under analysis.

\begin{equation}\label{eq:bond_stiffness_constant}
    c = 
    \begin{cases}
        \dfrac{12E}{\pi \delta^4} & \text{3D} \\[3ex]
        \dfrac{9E}{\pi t\delta^3} & \text{Plane stress} \\[3ex]
        \dfrac{48E}{5\pi t\delta^3} & \text{Plane strain} \\
    \end{cases}
\end{equation}

The evolution of the non-linear bond softening parameter $d$ is defined by Eq. (\ref{eq:exp_damage_parameter}). This function describes an exponentially decaying curve with a linear term. As the bond stretch $s$ approaches the critical stretch $s_c$, the linear term forces the softening curve to decay linearly and intersect with $s_c$. $\alpha$ controls the position of the transition from exponential to linear decay, and $k$ controls the rate of exponential decay. The linear elastic limit $s_0$ is defined empirically as $f_t/E$. This definition of $s_0$ is not objective but it has been shown to provide good results. 

\begin{equation}\label{eq:exp_damage_parameter}
    d = 1 - \frac{s_0}{s}\left(1 - \frac{1 - exp\left(-k\frac{s - s_0}{s_c - s_0}\right)}{1 - exp(-k)} + \alpha \left(1 - \frac{s - s_0}{s_c - s_0}\right) \right)(1 + \alpha)^{-1}
\end{equation}

The energy required to break a bond is defined by Eq. (\ref{eq:exp_bond_energy}). Only the energy consumed during the softening stage is considered (between the limits $s_0$ and $s_c$). It is important that the softening curve intersects with $s_c$ so that the integral in Eq. (\ref{eq:exp_bond_energy}) can be evaluated.

\begin{equation}\label{eq:exp_bond_energy}
    w = \int_{s_0}^{s_c} f \xi ds = \frac{cs_0 \xi (s_0 - s_c)(2k - 2e^k + \alpha k - \alpha ke^k + 2)}{2k(e^k - 1)(1 + \alpha)}
\end{equation}

\subsection*{Three-dimensional case}

Using Eq. (\ref{eq:exp_bond_energy}), the material fracture energy $G_F$ for a three-dimensional peridynamic body can be derived as follows

\begin{equation}
    G_F = \int_{0}^{\delta}\int_{0}^{2\pi}\int_{z}^{\delta}\int_{0}^{cos^{-1}(z/\xi)} w \xi^2 sin\phi \; d\phi \; d\xi \; d\theta \; dz
\end{equation}

\begin{equation}\label{eq:exp_fracture_energy}
    G_F = \frac{\pi c s_0 \delta^5 (s_0 - s_c)(2k - 2e^k + \alpha k - \alpha k e^k + 2)}{10k(e^k - 1)(1 + \alpha)}
\end{equation}

Rearranging Eq. (\ref{eq:exp_fracture_energy}), the critical stretch $s_c$ can be defined in terms of $s_0$, $k$ and $\alpha$. 

\begin{equation}\label{eq:exp_sc}
    s_c = \frac{10k(1 - e^k)(1 + \alpha)\left(G_F - \frac{\pi c \delta^5 s_0^2 (2k - 2e^k + \alpha k - \alpha k e^k + 2)}{10k(e^k - 1)}\right)}{c \delta^5 s_0 \pi (2k - 2e^k + \alpha k - \alpha k e^k + 2)}
\end{equation}

\subsection*{Two-dimensional case}

Using Eq. (\ref{eq:exp_bond_energy}), the material fracture energy $G_F$ for the two-dimensional case can be derived as follows

\begin{equation}\label{eq:exp_fracture_energy_2D}
    G_F = 2h\int_{0}^{\delta}\int_{z}^{\delta}\int_{0}^{cos^{-1}(z/\xi)} w \xi \; d\phi \; d\xi \; dz
\end{equation}

\noindent where $h$ is the thickness. Rearranging Eq. (\ref{eq:exp_fracture_energy_2D}), the critical stretch $s_c$ can be defined in terms of $s_0$, $k$ and $\alpha$. 

\begin{equation}\label{eq:exp_sc_2D}
    s_c= \frac{4k(1 - e^k)(1 + \alpha)\left(G_F - \frac{h c \delta^4 s_0^2 (2k - 2e^k + \alpha k - \alpha k e^k + 2)}{4k(e^k - 1)(1 + \alpha)}\right)}{h c \delta^4 s_0 (2k - 2e^k + \alpha k - \alpha k e^k + 2)}
\end{equation}

The proposed model provides an explicit definition of the critical stretch $s_c$ and an unambiguous relationship between $s_c$, $k$ and $\alpha$. 


\bibliographystyle{elsarticle-num-names}
\bibliography{references.bib} 

\end{document}